\documentclass[paper=a4, fontsize=11pt]{scrartcl}	
\usepackage{marvosym}
\usepackage{amsmath,amsfonts,amsthm,mathtools,pgfplots}
\usepackage{tikz-cd}
\usepackage[all]{xy}			
\usepackage{graphicx}
\usepackage{caption,subcaption}
\newtheorem{theorem}{Theorem}[section]
\newtheorem{lemma}[theorem]{Lemma}
\newtheorem{proposition}[theorem]{Proposition}

\newtheorem{remark}[theorem]{Remark}
\newtheorem{example}[theorem]{Example}
\newtheorem{definition}[theorem]{Definition}	
\usepackage{tikz}
\usepackage{mathdots}
\usepackage{enumerate}
\usepackage{xcolor}
\usepackage{adjustbox}
\tikzstyle{place}=[circle,draw=blue!50,fill=blue!20,thick]
\title{$\tilde{A}$ and $\tilde{D}$ type cluster algebras: Triangulated surfaces and friezes}
\author{Joe Pallister\footnote{JSPS International Research Fellow}\\\\Department of Mathematics and Informatics, Faculty of Science\\ Chiba University, Japan}
\begin{document}
\maketitle
\begin{abstract}
By viewing $\tilde{A}$ and $\tilde{D}$ type cluster algebras as triangulated surfaces, we find all cluster variables in terms of either (i) the frieze pattern (or bipartite belt) or (ii) the periodic quantities previously found for the cluster map associated with these frieze patterns. We show that these cluster variables form friezes which are precisely the ones found in \cite{assemdupont} by applying the cluster character to the associated cluster category.
\end{abstract}
\section{Introduction}
Cluster algebras were first defined in \cite{clusteri}, they are generated by cluster variables. Starting with an initial seed; a quiver and a set of cluster variables, there is an operation called mutation, at any vertex, which gives a new quiver and a new cluster variable. The set of cluster variables is given by applying all possible sequences of mutations.

By restricting so that we are only allowed to mutate at sinks or sources (where mutation becomes much simpler) we find cluster variables $X^i_n$ satisfying the frieze pattern:
\begin{equation}\label{friezeformulaintro}
X^i_nX^i_{n+1}=1+\left(\prod_{j\rightarrow i}(X^j_n)^{|b_{ji}|}\right)\left(\prod_{i\rightarrow j}(X^j_{n+1})^{|b_{ji}|}\right)
\end{equation}
for each vertex $i$ and $n\in\mathbb{Z}$. Here $b_{ji}=-b_{ij}$ is the number of arrows from $j$ to $i$ in $Q$. We consider (\ref{friezeformulaintro}) as an equation for determining $X^i_{n+1}$. Remarkably, if the initial quiver is of Dynkin type, then this isn't a restriction at all, as proved in \cite{clustersiv}.
\begin{theorem}
If the quiver $Q$ is of Dynkin type then there are finitely many cluster variables and each of them can be found on the frieze pattern.
\end{theorem}
If $Q$ if not Dynkin then the frieze pattern does not contain all of the cluster variables. The aim of this work is to describe what appears outside of the frieze pattern for the affine quivers $\tilde{A}$ and $\tilde{D}$. 

The frieze pattern was considered in \cite{kellerscherotzke} where, via the cluster category, the authors proved that linear relations exist between the frieze variables if and only if $Q$ is Dynkin or affine type. A more direct construction of these linear relations was given in \cite{fordyhone,pallisterlinear} by proving that there exists periodic quantities $J_n$ and $\tilde{J}_n$ for $\tilde{A}$ type and $J'_n$ for $\tilde{D}$ type. These are functions of the $X^i_n$ that are fixed by sending $n \mapsto n+a$ for an appropriate $a$. For example, in $\tilde{D}$ case, we have
\[
J'_n:=\frac{X^1_{n+1}+X^1_{n-1}}{X^2_n}=\frac{X^1_{n+N-1}+X^1_{n+N-3}}{X^2_{n+N-2}}=J'_{n+N-2}
\]
where the vertices $1$ and $2$, that appear as superscripts, are given in Figure \ref{Dquiver}. The first result of this paper is that these periodic quantities are also cluster variables, but do not live on the frieze pattern. 

In order to find the rest of the cluster variables we rely on a method pioneered in \cite{clustersandtriangulated1} for viewing cluster variables as arcs of triangulations of particular surfaces. For $\tilde{A}_{q,p}$ the surface is an annulus with $q$ and $p$ marked points on either boundary component. An example is shown in Figure \ref{exampletriangulation} which gives an annulus by gluing along the dotted lines. For $\tilde{D}_N$ the surface is a disc with two internal marked points and $N-2$ marked points on the boundary, shown in Figure \ref{Dtypebefore}. 

By finding all of the possible arcs for $\tilde{A}$ and $\tilde{D}$ types, we prove that all cluster variables can be described in terms of the frieze sequence (\ref{friezeformulaintro}) and determinants:
\begin{equation}\label{determinantintro}
D^m_a(F_n):=
\begin{vmatrix}
F_n & 1 & 0 \\
1 & F_{n+a} & 1 & 0 \\
0 & 1 & F_{n+2a} & 1 & 0 \\
& 0 & 1 & F_{n+3a} & \ddots \\
& & 0 & \ddots & \ddots & 1 \\
& & & & 1 & F_{n+(m-1)a}
\end{vmatrix}
\end{equation}
where $F$ is a function of the periodic quantities as described in the following two theorems.
\begin{theorem}\label{Atypearcs}
\begin{enumerate}[(i)]
For $\tilde{A}_{q,p}$ cluster algebras the arcs on the annulus are of three types:
\item
The arcs that connect the two boundary components. These are in bijection with the frieze variables $X^i_n$. 
\item
The arcs connecting the boundary component with $q$ marked points to itself, in bijection with the cluster variables $D^l_p(J_{jp})$, for $j=0,\ldots,q-1$ and $l=1,\ldots,q-1$.
\item 
The arcs connecting the other boundary component (with $p$ marked points) to itself, in bijection with $D^l_q(\tilde{J}_{jq})$, for $j=0,\ldots,p-1$ and $l=1,\ldots,p-1$.
\end{enumerate}
\end{theorem}
The corresponding result for $\tilde{D}$ type is:
\begin{theorem}
\begin{enumerate}[(i)]
For $\tilde{D}_{N}$ cluster algebras we divide the arcs of the twice punctured disk into four types.
\item
The arcs that connect the boundary vertices such that the two punctures lie on different sides of this arc.  
\item
The arcs connecting the boundary to the punctures. These and the arcs of (i) are in bijection with the frieze variables $X^i_n$.
\item 
The arcs connecting the boundary component to itself, in bijection with the $D^l_1({J'}_{j})$, for $j=0,\ldots,N-3$ and $l=1,\ldots, N-3$.
\item
The three exceptional arcs $\Gamma_{\mathrm{except}}$:
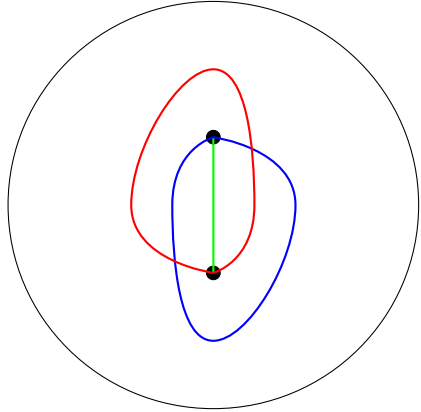
\begin{equation}\label{exceptionalarcsintro}
\begin{tikzpicture}[scale=0.9, every node/.style={fill=white}]
\draw (0,0) circle [radius=3.0];
\draw[fill=black] (0,1) circle [radius=0.1cm];
\draw[fill=black] (0,-1) circle [radius=0.1cm];
\draw[green, thick] plot [smooth, tension=1] coordinates { (0,1) (0,-1) };
\draw[blue, thick] plot [smooth, tension=1] coordinates { (0,1) (-0.6,0) (0,-2) (1.2,0) (0,1) };
\draw[red, thick] plot [smooth, tension=1] coordinates {(0,-1) (-1.2,0) (0,2) (0.6,0) (0,-1)};
\end{tikzpicture}
\end{equation}
\end{enumerate}
\end{theorem}
If we forget about arcs and surfaces these two theorems give a description of all of the cluster variables in $\tilde{A}$ and $\tilde{D}$ type: 
\begin{theorem}\label{findingallclustervarstheorem}
For $\tilde{A}_{q,p}$ cluster algebras the cluster variables are
\begin{equation}
\left\{X^i_n\:\middle|\:
\begin{aligned}
&i=1,\ldots,N \\ &n\in\mathbb{Z}
\end{aligned}
\right\}
\cup
\left\{D^l_p(J_{jp})\:\middle|\:
\begin{aligned}
j=0,\ldots,q-1 \\ l=1,\ldots,q-1
\end{aligned}
\right\}
\cup
\left\{D^l_q(\tilde{J}_{jq})\:\middle|\:
\begin{aligned}
j=0,\ldots,p-1 \\ l=1,\ldots,p-1
\end{aligned}
\right\}
\end{equation}
For $\tilde{D}_N$ the cluster variables are 
\begin{equation}
\left\{X^i_n\:\middle|\:
\begin{aligned}
&i=1,\ldots,N+1 \\ &n\in\mathbb{Z}
\end{aligned}
\right\}
\cup
\left\{D^l_1(J'_{j})\:\middle|\:
\begin{aligned}
j=0,\ldots,N-3 \\ l=1,\ldots,N-3
\end{aligned}
\right\}
\cup
\Gamma_{\mathrm{except}}
\end{equation} 
\end{theorem}
Where $\Gamma_{\mathrm{except}}$ is the set of cluster variables associated with three exceptional arcs (\ref{exceptionalarcsintro}).

Friezes were defined in \cite{coxeterfrieze} where they were shown to have connections with continued fractions and Farey series. Further links with diverse topics including triangulations of polygons \cite{conway1973triangulated}, Auslander-Reiten theory \cite{gabriel1980auslander} and moduli spaces of points in projective space \cite{morier20122} were later found. They were shown to be linked with cluster algebras \cite{calderochapoton} via the cluster category and a more direct link was studied in \cite{frises}. Here we prove that in $\tilde{A}$ and $\tilde{D}$ cases the off-frieze pattern cluster variables form (slightly generalised) friezes.
\begin{proposition}\label{friezepropintro}
In $\tilde{A}_{q,p}$ and $\tilde{D}_N$ cluster algebras the non-frieze pattern cluster variables defined by (\ref{determinantintro}) form friezes:
\begin{equation}\setcounter{MaxMatrixCols}{20}
\begin{matrix}
\ldots & & 1 & & 1 & & 1 & &  \\
& D^1_a(F_n) & & D^1_a(F_{n+a}) & & D^1_a(F_{n+2a}) & & \ldots & & \\
\ldots & & \mathcal{D}^2_{a}(F_n) & & \mathcal{D}^2_{a}(F_{n+a}) & & \mathcal{D}^2_{a}(F_{n+2a}) & &  \\
& \mathcal{D}^3_{a}(F_{n-a}) & & \mathcal{D}^3_{a}(F_n) & & \mathcal{D}^3_{a}(F_{n+a}) & & \ldots \\
\ldots & & \mathcal{D}^4_{a}(F_{n-a}) & & \mathcal{D}^4_{a}(F_{n}) & & \mathcal{D}^4_{a}(F_{n+a}) & &  \\ \\
 & & \vdots & & \vdots && \vdots  \\ \\
& \mathcal{D}^L_{a}(F_{n-3a}) & & \mathcal{D}^L_{a}(F_{n-2a}) & & \mathcal{D}^L_{a}(F_{n-a}) & & \ldots
\end{matrix}
\end{equation}
i.e. each diamond 
$\begin{matrix}
& \beta & \\
\alpha & & \delta \\
& \gamma &
\end{matrix}$
satisfies $\alpha\delta-\beta\gamma=1$. In $\tilde{A}_{q,p}$ type there are two of these friezes, given by
\[
(i) \qquad F_n=J_{jp}, \quad a=p, \quad L=q-1.
\]
\[
(ii) \qquad F_n=\tilde{J}_{jq}, \quad a=q, \quad L=p-1.
\]
In $\tilde{D}_N$ type the single frieze is given by 
\[
F_n=J'_{j}, \quad a=1, \quad L=N-3.
\]
\end{proposition}
For these quivers this structure has previously been found \cite{assemdupont} in terms of the cluster category. There is a map $X_{?}$ from the objects of the cluster category to the cluster algebra, such that the frieze pattern (\ref{friezeformulaintro}) is the image of the transjective component. It was also shown that the exceptional tubes form friezes under this map, and we link our work to this construction.
\begin{proposition}\label{intropropclusterfrieze}
In $\tilde{A}_{q,p}$ type the frieze pattern (\ref{friezeformulaintro}) and the friezes of Proposition \ref{friezepropintro} give precisely the cluster frieze given in \cite{assemdupont}. In $\tilde{D}_N$ type the frieze pattern (\ref{friezeformulaintro}) and Proposition \ref{friezepropintro} give the cluster frieze except for the portion on the period $2$ tubes, which we have not constructed.
\end{proposition}
In this way we can see the $\tilde{A}$ and $\tilde{D}$ type cluster algebras as a union of friezes ($3$ friezes in the $\tilde{A}$ case and $2$ in the $\tilde{D}$ case) linked by the relation between the expression in the frieze variables defining the periodic quantities $J$, $\tilde{J}$ and $J'$.

This paper is arranged as follows:

In Section \ref{backgroundsection} we explain the background material we need, including the construction of the frieze pattern (\ref{friezeformulaintro}) and some properties of the periodic quantities we will use to prove the determinant formula (\ref{determinantintro}). We then describe how certain cluster algebras can be viewed as surfaces and how mutation works in this picture. Next we show that the $X^i_n$ defined by (\ref{friezeformulaintro}) form friezes on repetition quivers. We then define the cluster category and the cluster character and show that for $\tilde{A}$ and $\tilde{D}$ quivers the frieze (\ref{friezeformulaintro}) agrees with the frieze constructed on the transjective component of the cluster category via the cluster character, as given in \cite{assemdupont}. Finally we give the definition of a cluster frieze, which is a frieze given on the whole cluster category, not just the transjective component.

In Section \ref{Friezeconstructionsection} we construct the friezes of Proposition \ref{friezepropintro}, while not yet proving that the frieze entries are cluster variables.

In Section \ref{Atypesection} we describe the possible arcs (cluster variables) in $\tilde{A}$ type, proving the $\tilde{A}$ parts of Theorems \ref{findingallclustervarstheorem} and \ref{Atypearcs} and Proposition \ref{friezepropintro}. We then compare our friezes with those constructed in \cite{assemdupont}, proving Proposition \ref{intropropclusterfrieze} in the $\tilde{A}$ case. 

Section \ref{Dtypesection} is analogous to Section \ref{Atypesection}, but for $\tilde{D}$ type.
\section{Review of cluster algebras, the cluster map, cluster algebras as triangulations of surfaces, friezes and the cluster category}\label{backgroundsection}
Here we describe the disparate elements that combine to give our results. Firstly we give the construction of cluster mutation and cluster algebras. We then discuss how the cluster map is defined for period $1$ quivers and  bipartite quivers before giving a general definition that gives the frieze pattern (\ref{friezeformulaintro}).

Next we give previous results that have appeared for the frieze pattern in $\tilde{A}$ and $\tilde{D}$ types, including identifying periodic quantities for their respective cluster maps. We then show how certain cluster algebras can be described as triangulations of surfaces, since we aim to interpret these periodic quantities as arcs in these triangulations.

We next define friezes and show how they relate to cluster algebras.

Finally we give a very brief introduction to the cluster category and the cluster character, a map from the objects of this category to the set of cluster variables. We describe the cluster categories in $\tilde{A}$ and $\tilde{D}$ type and describe the frieze constructed in \cite{assemdupont} on the objects.
\subsection{Cluster algebras}
Here we give our definition of cluster algebras. In this paper we mean a cluster algebra without coefficients or frozen variables. 

A quiver $Q$ is a directed graph where multiple edges are allowed. We disallow loops or 2-cycles. Quiver mutation $\mu_k$ at any vertex $k$ is defined in $3$ steps: 
\begin{enumerate}
\item For each length two path $i\rightarrow k \rightarrow j$ add a new arrow $i\rightarrow j$.
\item Reverse the direction of all arrows entering or exiting $k$.
\item Delete all 2-cycles that have appeared.
\end{enumerate}
This gives a new quiver $\mu_k(Q)$. 

The exchange matrix for a quiver $Q$ is the skew-symmetric matrix $B$ with entries $b_{ij}$, the number of arrows from $i$ to $j$, with $b_{ji}=-b_{ij}$. The mutation $\mu_k$ acts on the $b_{ij}$ as
\begin{equation}\label{mutatedBmatrixentries}
\mu_k(b_{ij})=
\begin{cases}
-b_{ij} & \mbox{ if } i=k \mbox{ or } j=k, \\
b_{ij}+\frac{1}{2}(|b_{ik}|b_{kj}+b_{ik}|b_{kj}|) & \mbox{ otherwise. }
\end{cases}
\end{equation}
In addition to this, we will also attach a cluster variable $x_i$ at each vertex $i$. Cluster mutation at vertex $k$, also denoted $\mu_k$, fixes all variables $x_i$ with $i\neq k$ but transforms $x_k$ as follows:
\begin{equation}\label{clustermutationformula}
\mu_k(x_k):=\frac{1}{x_k}\left(\prod_{i\rightarrow k} x_i+\prod_{i\leftarrow k}x_i\right).
\end{equation}
Here the two products run over the arrows into and out of $k$ respectively, e.g. in the first product we have an $x_i$ for every arrow from $i$ to $k$. We consider $\mu_k$ as a mutation both of the quiver and of the cluster variables.

The cluster algebra $\mathcal{A}_Q$ is the algebra over $\mathbb{Z}$ generated by the cluster variables obtained by any sequence of mutations applied to $Q$.

\subsection{Periodic quivers and the cluster map}\label{dynamicalsystemssubsection}
In this section we define periodic quivers and (a restricted version of) the cluster map that follows. We show how the cluster map has been extended to bipartite quivers before giving a general definition for all acyclic quivers. We then see how the frieze pattern (\ref{friezeformulaintro}) follows from this map.

Period $1$ quivers were defined and classified in \cite{fordymarsh}. These are $N$ vertex quivers $Q$, labelled $0,1,\ldots,N-1$, satisfying $\mu_0(Q)=\rho(Q)$, where $\rho=(0,N,N-1,\ldots,2,1)$ is a permutation acting on the vertices of $Q$. This means that mutation at vertex $0$ is tantamount to a (specific) relabelling of the vertices of $Q$. 

By taking initial cluster variables $x_0,x_1,\ldots,x_{N-1}$ mutation at $0$ will give a new cluster variable, which we call $x_N$, determined by
\begin{equation}\label{firstperiod1mutation}
x_Nx_0=F(x_{N-1},x_{N-2},\ldots,x_2,x_1)
\end{equation}
for an appropriate function $F$. Next mutation at $1$ will give a new cluster variable called $x_{N+1}$ satisfying
\[
x_{N+1}x_1=F(x_N,x_{N-1},\ldots,x_2)
\] 
where, due to the quiver being period $1$, this $F$ is the same as in (\ref{firstperiod1mutation}). Continuing this process gives the recurrence
\[
x_{n+N}x_n=F(x_{n+N-1},x_{n+N-2},\ldots,x_{n+2},x_{n+1})
\] 
for all $n$. An example of this is given by an orientation of an affine $A$ type diagram with $p$ and $q$ arrows pointing clockwise and anticlockwise, respectively, called $\tilde{A}_{q,p}$. This is given by taking the diagram shown in Figure \ref{Atypediagram}, reducing the labels modulo $N$, and orienting so that each arrow points from the lower label to the higher.
\begin{lemma}\label{Atypeconstructionlemma}
The $\tilde{A}_{q,p}$ quiver, with $p$ and $q$ coprime, gives the recurrence 
\begin{equation}\label{Atyperecurrence}
x_{n+p+q}x_n=x_{n+q}x_{n+p}+1.
\end{equation}
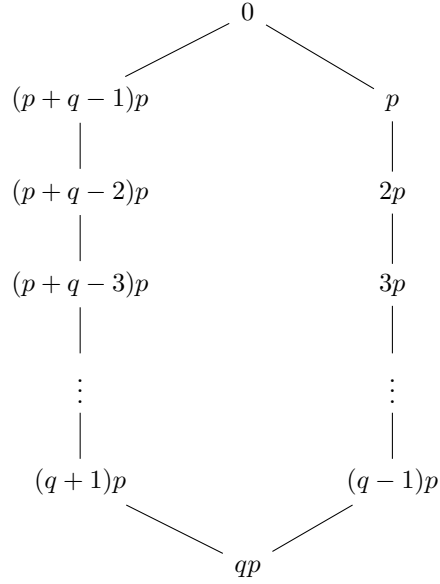
\begin{figure}
\centering
\begin{equation*}
\adjustbox{scale=0.9,center}{%
\begin{tikzcd}
& 0\arrow[dr, dash] \arrow[dl, dash] & \\
(p+q-1)p\arrow[d,dash] &  & p\arrow[d,dash] \\
(p+q-2)p\arrow[d,dash] &  & 2p\arrow[d,dash] \\
(p+q-3)p\arrow[d,dash] & & 3p\arrow[d,dash] \\
\vdots\arrow[d,dash] & & \vdots\arrow[d,dash] \\
(q+1)p\arrow[dr,dash] & & (q-1)p\arrow[dl,dash] \\
& qp &
\end{tikzcd}
}
\end{equation*}
\caption{The $\tilde{A}$ type diagram used to obtain the $\tilde{A}$ type recurrence relation.}
\label{Atypediagram}
\end{figure}
\end{lemma}
\begin{example}
For $p=7$ and $q=8$ the recurrence of this lemma is given by the quiver shown in Figure \ref{EuclideanA15}.
\begin{figure}[h]
\centering
\begin{equation*}
\begin{tikzcd}
9 & 1\arrow[l]\arrow[r] & 8 & 0\arrow[l]\arrow[r] & 7\arrow[r] & 14 & \arrow[l]6\arrow[d] \\
2\arrow[u]\arrow[d] &   &   &   &   &    & 13 \\
10 & 3\arrow[l]\arrow[rr] && 11 & 4\arrow[l]\arrow[r] & 12 & 5\arrow[l]\arrow[u]
\end{tikzcd}
\end{equation*}\caption{Our orientation of the $\tilde{A}_{8,7}$ quiver.} \label{EuclideanA15}
\end{figure}
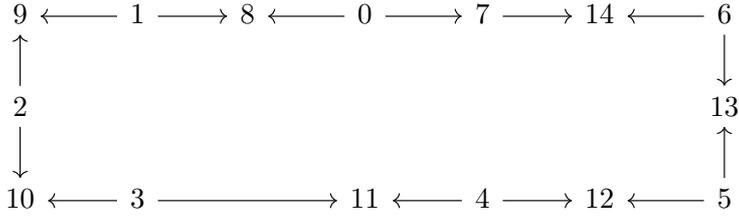
\end{example}

The map defined by these mutations
\[
\varphi:(x_n,x_{n+1},\ldots,x_{n+N-1})\mapsto (x_{n+1},x_{n+2},\ldots,x_{n+N})
\]
is know as the cluster map. In \cite{pallisterlinear} it was generalised to include bipartite quivers $Q$; we take $\mu_{\mathrm{sink}}$ and $\mu_{\mathrm{source}}$ to be the compositions of mutations at every sink and source in $Q$, respectively, then the product $\mu:=\mu_{\mathrm{source}}\circ\mu_{\mathrm{sink}}$ fixes the quiver but gives new cluster variables. We let $X^i_n$ be the cluster variable at vertex $i$ after $n$ applications of $\mu$ (with $n=0$ giving the initial variables) then the map
\begin{equation}\label{generalisedclustermap}
\varphi:(X^1_n,X^2_n,\ldots,X^{N+1}_n)\mapsto (X^1_{n+1},X^2_{n+1},\ldots,X^{N+1}_{n+1})
\end{equation}
is also known as a cluster map.

For a general acyclic quiver this map is obtained by first taking $X^i_0$ for each vertex $i$. We then take $X^i_n$, for $n\neq 0$, to satisfy
\begin{equation}\label{friezeformula}
X^i_nX^i_{n+1}=1+\left(\prod_{j\rightarrow i}(X^j_n)^{|b_{ji}|}\right)\left(\prod_{i\rightarrow j}(X^j_{n+1})^{|b_{ji}|}\right)
\end{equation}
where the products is taken over all arrows in $Q$. This is called a ``generalised frieze pattern" in \cite{kellerscherotzke}. It can be shown that all of the $X^i_n$ can be obtained by cluster mutation, similarly to the construction of (\ref{generalisedclustermap}). The $\tilde{A}$ type cluster variables obtained by (\ref{Atyperecurrence}) are also included here, with the identification $X^i_n\leftrightarrow x_{i-nN}$. 
\subsection{Periodic quantities for the cluster map}
Here we examine the periodic quantities found for the cluster map where the quiver is of affine type. These immediately give linear relations between the frieze variables with periodic coefficients and, with some work, linear relations with constant coefficients can be obtained. Here we explain how this is done.

In $\tilde{A}$ type, both papers \cite{fordymarsh,fordyhone} prove that the recurrence (\ref{Atyperecurrence}) has the periodic quantities 
\begin{equation}\label{periodicquantities}
J_n:=\frac{x_{n+2p}+x_n}{x_{n+p}}, \qquad \tilde{J}_n:=\frac{x_{n+2q}+x_n}{x_{n+q}}
\end{equation}
with period $q$ and $p$, respectively. By this we mean that $J_{n+q}=J_n$ and $\tilde{J}_{n+p}=\tilde{J}_n$. 
Immediately we see that the $x_n$ satisfy the linear relations
\begin{equation}\label{periodiclinearrelations}
x_{n+2p}-J_nx_{n+p}+x_n=0, \qquad x_{n+2q}-\tilde{J}_nx_{n+q}+x_n=0
\end{equation}
with periodic coefficients. The authors of \cite{fordyhone} then use this to construct linear relations with constant coefficients.
\begin{theorem}
The cluster variables $x_n$ satisfy the linear relation 
\[
x_{n+2qp}-\mathcal{K}x_{n+qp}+x_n=0
\]
where $\mathcal{K}$ is constant.
\begin{proof}
We define the matrices
\[
\Psi_n:=
\begin{pmatrix}
x_{n+p+q} & x_{n+q} \\
x_{n+p} & x_{n}
\end{pmatrix}
\qquad 
L_n:=
\begin{pmatrix}
J_n & 1 \\
-1 & 0 
\end{pmatrix}
\qquad 
\tilde{L}_n:=
\begin{pmatrix}
\tilde{J}_n & -1 \\
1 & 0
\end{pmatrix}
\]
such that, by (\ref{periodiclinearrelations}), $\Psi_{n+p}=\Psi_nL_n$ and $\Psi_{n+q}=\tilde{L}_n\Psi_n$. By defining
\[
M_n:=L_nL_{n+p}L_{n+2p}\ldots L_{n+(q-1)p}, \qquad \tilde{M}_n:=\tilde{L}_{n+(p-1)q}\ldots\tilde{L}_{n+2q}\tilde{L}_{n+q}\tilde{L}_n
\]
we have $\Psi_nM_n=\Psi_{n+qp}$ and $\tilde{M}_n\Psi_n=\Psi_{n+qp}$. We then apply the Cayley-Hamilton theorem to $M_n$:
\[
M_n^2-\mathcal{K}M_n+I=0
\]
where ${\mathcal{K}=\mathrm{Tr}(M_n)=\mathrm{Tr}(\tilde{M}_n)}$ is constant. Multiplying this matrix equation by $\Psi_n$ on the left gives the linear relation.
\end{proof}
\end{theorem} 
The first question we ask is about the entries of $M_n$ and $\tilde{M}_n$. We define the matrices $M^m_n$ as the product of the first $m$ of these $L$ matrices:
\begin{equation}\label{matrixproduct}
M^m_n:=L_nL_{n+p}L_{n+2p}\ldots L_{n+(m-1)p}.
\end{equation}
For example $M^3_n$ is given by
\[
M^3_n=
L_nL_{n+p}L_{n+2p}=
\begin{pmatrix}
J_{n+2p}J_{n+p}J_n-J_{n+2p}-J_n & J_{n+p}J_n-1 \\
-J_{n+2p}J_{n+p}+1 & -J_{n+p}
\end{pmatrix}
\]
Similarly we let
\[
\tilde{M}^m_n:=\tilde{L}_{n+(m-1)q}\ldots\tilde{L}_{n+2q}\tilde{L}_{n+q}\tilde{L}_n.
\]
We prove the following in Section \ref{Friezeconstructionsection}:
\begin{proposition}\label{matrixentryclustervariables}
The matrices $M^m_n$ are given by
\[
M^m_n=
\begin{pmatrix}
A^m_n & A^{m-1}_n \\
-A^{m-1}_{n+p} & -A^{m-2}_{n+p}
\end{pmatrix}
\]
where $A^m_n$ satisfies the recurrence in $m$:
\begin{equation}\label{mainreccurence}
A^m_n=J_{n+(m-1)p}A^{m-1}_n-A^{m-2}_n, \qquad A^1_n=J_n, \qquad A^0_n=1
\end{equation}
for each $n$. Conversely the matrices $\tilde{M}^m_n$ are given by 
\[
\tilde{M}^m_n=
\begin{pmatrix}
\tilde{A}^m_n & -\tilde{A}^{m-1}_{n+q} \\
\tilde{A}^{m-1}_n & -\tilde{A}^{m-2}_{n+q}
\end{pmatrix}
\]
satisfying
\[
\tilde{A}^m_n=\tilde{J}_{n+(m-1)q}\tilde{A}^{m-1}_n-\tilde{A}^{m-2}_n, \qquad \tilde{A}^1_n=\tilde{J}_n, \qquad \tilde{A}^0_n=1.
\]
\end{proposition}
In \cite{pallisterlinear} periodic quantities were found for the cluster map (\ref{generalisedclustermap}) for $\tilde{D}$ and $\tilde{E}$ quivers, as shown in Figure \ref{periodictable}. \begin{figure}
\begin{center}
\begin{tabular}{|c|c|c|c|}
\hline 
Quiver & Period  & Explicit Expression in \cite{pallisterlinear} \\
\hline\noalign{\smallskip}
& $N-2$ & Equation 22\\
$\tilde{D}_N$ & $2$ & Lemma 4.1 \\
& $2$  & Lemma 4.1 \\
\hline\noalign{\smallskip}
& $3$ & Lemma 4.9 \\
$\tilde{E}_6$ & $3$  & Lemma 4.9  \\
& $2$ &  Lemma 4.12 \\
\hline\noalign{\smallskip}
& $4$  & Theorem 4.20 \\
$\tilde{E}_7$ & $3$ & Theorem 4.19 \\
& $2$ & Equation 42 \\
\hline\noalign{\smallskip}
 & $5$ & Lemma 4.26 \\
$\tilde{E}_8$ & $3$ & Theorem 4.28 \\
 & $2?$  & Conjecture 4.30 \\
\hline
\end{tabular}
\end{center}\caption{Periodic quantities found for the cluster map for $\tilde{D}$ and $\tilde{E}$ quivers.}\label{periodictable}
\end{figure}
The linear relations that follow from these were also obtained in \cite{kellerscherotzke} using the cluster category. For $\tilde{D}$ type, the period $N-2$ quantity can be written
\begin{equation}\label{DtypeJintro}
J'_n:=\frac{X^1_{n+1}+X^1_{n-1}}{X^2_n}
\end{equation}
and can also be used to give linear relations with constant coefficients, as in the $\tilde{A}$ case.
\begin{proposition}
In the $\tilde{D}_N$ cases, for even $N$ the cluster variables at vertex $1$ satisfy the constant coefficient linear relation:
\[
X^1_{n+2N-4}-\mathcal{K}'X^1_{n+N-2}+X^1_n=0
\]
and for odd $N$:
\[
X^1_{n+4N-8}-\mathcal{K}'X^1_{n+2N-4}+X^1_n=0
\]
where vertex $1$ is as shown in Figure \ref{Dquiver}.
\end{proposition}
The proof of this is similar to the $\tilde{A}$ case. For $N$ even the constant $\mathcal{K}'$ is given by the trace of
\begin{equation}\label{Deven}
{L'}_n{L'}_{n+1}\ldots{L'}_{n+N-3}, \qquad 
{L'}_n:=
\begin{pmatrix}
0 & -1 \\
1 & J'_n
\end{pmatrix}
\end{equation}
and for odd $N$ the trace of 
\begin{equation}\label{Dodd}
{L'}_n{L'}_{n+1}\ldots{L'}_{n+2N-5}.
\end{equation}
For the same ${L'}_n$. Interestingly, the same product of matrices appears as in the $\tilde{A}$ case. A similar result holds too.
\begin{proposition}\label{Dtypematrixentries}
The entries of the matrix products (\ref{Deven}) and (\ref{Dodd}) are determined by their upper left entry, which we call ${A'}^m_n$ in either case, which satisfies the linear recurrence 
\begin{equation}\label{Dtypematrixlinear}
{A'}^m_n=J'_{n+m-1}{A'}^{m-1}_n-{A'}^{m-2}_n, \qquad {A'}^1_n=J'_n, \qquad {A'}^0_n=1
\end{equation}
for each $n$.
\end{proposition}
The proof of this is similar to the one for $\tilde{A}$ type given in Section \ref{Friezeconstructionsection}. Due to this result we ask how the $A^m_n$, $\tilde{A}^m_n$ and ${A'}^m_n$ fit into their respective cluster algebras. Are they are cluster variables and, if so, what cluster variables remain outside of the $A^m_n$, $\tilde{A}^m_n$ and ${A'}^m_n$ and the variables obtained by the cluster map? To answer these we need an alternative way of viewing these cluster algebras.
\subsection{Cluster algebras as triangulated surfaces}
In \cite{clustersandtriangulated1} the authors describe how to view both Dynkin and affine $A$ and $D$ type cluster algebras as triangulations of surfaces with marked points. In this section we briefly explain how this correspondence works: we describe how a quiver is obtained from a triangulation and what the analogues of quiver and cluster mutation are.

\begin{figure}
\centering
\begin{tikzpicture}
\draw[fill=black] (0,0) circle [radius=0.1cm];
\draw[fill=black] (1,0) circle [radius=0.1cm];
\draw[fill=black] (2,0) circle [radius=0.1cm];
\node at (3,0) {$\ldots$};
\draw[fill=black] (4,0) circle [radius=0.1cm];
\draw[fill=black] (0.5,2) circle [radius=0.1cm];
\draw[fill=black] (1.5,2) circle [radius=0.1cm];
\node at (2.5,2) {$\ldots$};
\draw[fill=black] (3.5,2) circle [radius=0.1cm];
\foreach \x in {-1,...,1}
\draw [-] (\x,0) to (\x+1,0);
\draw [-] (2,0) to (2.5,0);
\draw [-] (3.5,0) to (4,0);
\foreach \x in {4,...,6}
\draw [-] (\x,0) to (\x+1,0);
\draw [-] (7,0) to (7.5,0);
\draw [-] (8.5,0) to (10,0);
\foreach \x in {-1,...,1}
\draw [-] (\x,2) to (\x+1,2);
\foreach \x in {3,...,6}
\draw [-] (\x,2) to (\x+1,2);
\node at (0,-0.3) {0};
\node at (1,-0.3) {1};
\node at (2,-0.3) {2};
\node at (4,-0.3) {$q-1$};
\node at (0.5,2.3) {$q$};
\node at (1.5,2.3) {$q+1$};
\node at (3.5,2.3) {$p+q-1$};
\draw[fill=black] (5,0) circle [radius=0.1cm];
\draw[fill=black] (6,0) circle [radius=0.1cm];
\draw[fill=black] (7,0) circle [radius=0.1cm];
\node at (8,0) {$\ldots$};
\draw[fill=black] (9,0) circle [radius=0.1cm];
\draw[fill=black] (5.5,2) circle [radius=0.1cm];
\draw[fill=black] (6.5,2) circle [radius=0.1cm];
\node at (7.5,2) {$\ldots$};
\draw[fill=black] (8.5,2) circle [radius=0.1cm];
\foreach \x in {7,...,8}
\draw [-] (10+\x,0) to (10+\x+1,0);
\draw [-] (8,2) to (10,2);
\node at (5,-0.3) {0};
\node at (6,-0.3) {1};
\node at (7,-0.3) {2};
\node at (9,-0.3) {$q-1$};
\node at (5.5,2.3) {$q$};
\node at (6.5,2.3) {$q+1$};
\node at (8.5,2.3) {$p+q-1$};
\draw[dotted] (-0.5,-1) to (-0.5,3);
\draw[dotted] (4.5,-1) to (4.5,3);
\draw[dotted] (9.5,-1) to (9.5,3);
\draw[fill=black] (10,0) circle [radius=0.1cm];
\node at (10,-0.3) {$1$};
\draw plot [smooth, tension=1] coordinates {(3.5,2) (10,0)};
\draw plot [smooth, tension=1] coordinates {(-1.5,2) (5,0)};
\draw[fill=black] (-1.5,2) circle [radius=0.1cm];
\node at (-1.5,2.3) {$p+q-1$};
\draw [-] (-2,2) to (0,2);
\draw [-] (10,0) to (10.5,0);
\draw [-] (8.5,2) to (10,1.53);
\draw plot [smooth, tension=1] coordinates {(0,0) (2,0.75) (4,0)};
\draw plot [smooth, tension=1] coordinates {(1,0) (2.5,0.5) (4,0)};
\draw plot [smooth, tension=1] coordinates {(6,0) (7.5,0.5) (9,0)};
\draw plot [smooth, tension=1] coordinates {(5,0) (7,0.75) (9,0)};
\draw plot [smooth, tension=1] coordinates {(0.5,2) (2,1.25) (3.5,2)};
\draw plot [smooth, tension=1] coordinates {(5.5,2) (7,1.25) (8.5,2)};
\end{tikzpicture}\caption{The annulus associated with the $\tilde{A}_{q,p}$ cluster algebra, with some possible arcs.}\label{annuluswithsomearcs}
\end{figure}
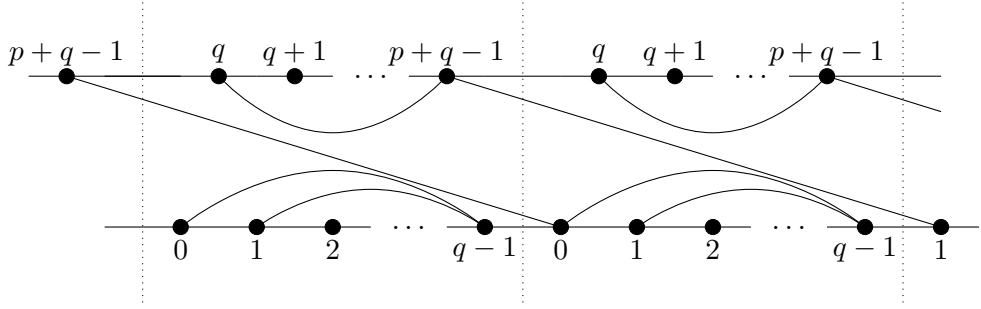

For $\tilde{A}_{q,p}$ cluster algebras we shall be dealing with triangulations of annuli with $q$ and $p$ marked points on opposite boundary components, an example of which can be seen in Figure \ref{exampletriangulation}, which gives an annulus by gluing along the dotted lines. In addition to the traditional triangles seen for this surface we also allow ``self-folded" triangles:
\begin{equation}\label{selffolded}
\scalebox{0.75}{
\begin{tikzpicture}[every node/.style={fill=white}]
\draw[fill=black] (0,0) circle [radius=0.15cm];
\draw[fill=gray] (0,2) circle [radius=0.15cm];
\draw [-] (0,0) to (0,2);
\draw plot [smooth, tension=1] coordinates {(0,0) (-0.9,2) (0,3.4) (0.9,2) (0,0)};
\node[scale=1.2] at (0,1) {$i$};
\node[scale=1.2] at (0,3.4) {$k$};
\end{tikzpicture}
}
\end{equation}
that will appear in the $\tilde{D}$ type triangulations, as in Figure \ref{Dtypebefore}. For $\tilde{A}$ type some possible arcs (not a full triangulation) are displayed in Figure \ref{annuluswithsomearcs}. 
To return the quiver from the triangulation:
\begin{enumerate}[(i)]
\item Attach a vertex to every arc.
\item For every pair of vertices $i$ and $j$ from step (i) that are part of the same non-self-folded triangle we draw an arrow $i\mapsto j$ if, while travelling clockwise around the triangle, $j$ comes directly after $i$. If no arcs lie on a boundary then the situation looks as follows:
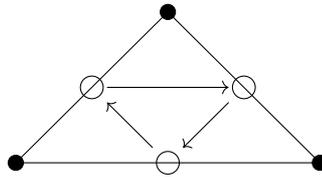
\[
\begin{tikzpicture}[every node/.style={fill=white}]
\draw[fill=black] (0,0) circle [radius=0.1cm];
\draw[fill=black] (4,0) circle [radius=0.1cm];
\draw[fill=black] (2,2) circle [radius=0.1cm];
\draw[fill=white] (1,1) circle [radius=0.15cm];
\draw[fill=white] (3,1) circle [radius=0.15cm];
\draw[fill=white] (2,0) circle [radius=0.15cm];
\draw [-] (0,0) to (2,2);
\draw [-] (0,0) to (4,0);
\draw [-] (2,2) to (4,0);
\draw [->] (1.2,1) to (2.8,1);
\draw [->] (1.8,0.2) to (1.2,0.8);
\draw [->] (2.8,0.8) to (2.2,0.2) ;
\end{tikzpicture} 
\]
where the white circles are the vertices of the quiver. If an arc lies on a boundary component then we do not draw a vertex there.
\item For every pair of arcs $k$ and $i$ forming a self-folded-triangle (\ref{selffolded}) and for every arrow $j\mapsto k$ we add an arrow $j\mapsto i$. Also for every arrow $k\mapsto j$ we add an arrow $i\mapsto j$.
\end{enumerate}
For arcs $k$ not part of self-folded triangles, quiver mutation is analogous to replacing the diagonal of a quadrilateral, $k$, with the other diagonal, $k'$:
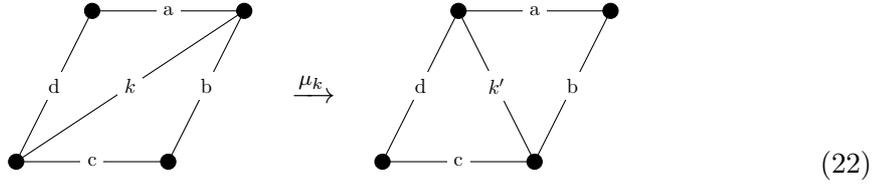
\begin{equation}\label{diagonalflip}
\begin{tikzpicture}[every node/.style={fill=white}]
\draw[fill=black] (0,0) circle [radius=0.1cm];
\draw[fill=black] (2,0) circle [radius=0.1cm];
\draw[fill=black] (1,2) circle [radius=0.1cm];
\draw[fill=black] (3,2) circle [radius=0.1cm];
\draw [-] (0,0) to (2,0);
\draw [-] (1,2) to (3,2);
\draw [-] (0,0) to node[scale=0.7] {$k$} (3,2);
\draw [-] (2,0) to (3,2);
\draw [-] (0,0) to (1,2);
\node[scale=0.7] at (2,2){a};
\node[scale=0.7] at (2.5,1){b};
\node[scale=0.7] at (1,0){c};
\node[scale=0.7] at (0.5,1){d};
\end{tikzpicture}
\quad
\begin{tikzpicture}[every node/.style={fill=white}]
\node at (0,0) {$\;$};
\node at (0,1) {$\xrightarrow{\mu_k}$};
\end{tikzpicture}
\quad
\begin{tikzpicture}[every node/.style={fill=white}]
\draw[fill=black] (0,0) circle [radius=0.1cm];
\draw[fill=black] (2,0) circle [radius=0.1cm];
\draw[fill=black] (1,2) circle [radius=0.1cm];
\draw[fill=black] (3,2) circle [radius=0.1cm];
\draw [-] (0,0) to (2,0);
\draw [-] (1,2) to (3,2);
\draw [-] (2,0) to node[scale=0.7] {$k'$} (1,2);
\draw [-] (2,0) to (3,2);
\draw [-] (0,0) to (1,2);
\node[scale=0.7] at (2,2){a};
\node[scale=0.7] at (2.5,1){b};
\node[scale=0.7] at (1,0){c};
\node[scale=0.7] at (0.5,1){d};
\end{tikzpicture}
\end{equation} 
known as a flip. Figures \ref{exampletriangulation} and \ref{triangulationaftermutation} show an example of this. In our work, for arcs $k$ that are the ``outside" of a self-folded triangle mutation looks like:
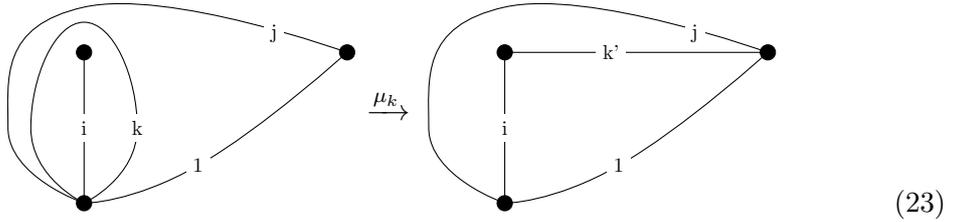
\begin{equation}\label{selffoldedmutation}
\begin{tikzpicture}[every node/.style={fill=white}]
\draw[fill=black] (0,-2) circle [radius=0.1cm];
\draw[fill=black] (0,-4) circle [radius=0.1cm];
\draw[fill=black] (3.46,-2) circle [radius=0.1cm];
\draw plot [smooth, tension=1] coordinates {(0,-4) (0,-2)};
\draw plot [smooth, tension=1] coordinates {(0,-4) (-0.7,-3) (0,-1.6) (0.7,-3) (0,-4)};
\draw plot [smooth, tension=0.9] coordinates {(0,-4) (-1,-3) (0,-1.4) (3.46,-2)};
\draw plot [smooth, tension=0.9] coordinates {(0,-4) (1.5,-3.5) (3.46,-2)};
\node[scale=0.7] at (0,-3){i};
\node[scale=0.7] at (0.7,-3){k};
\node at (4,-2.7){$\xrightarrow{\mu_k}$};
\node[scale=0.7] at (1.5,-3.5){1};
\node[scale=0.7] at (2.5,-1.75){j};
\end{tikzpicture} 
\begin{tikzpicture}[every node/.style={fill=white}]
\draw[fill=black] (0,-2) circle [radius=0.1cm];
\draw[fill=black] (0,-4) circle [radius=0.1cm];
\draw[fill=black] (3.46,-2) circle [radius=0.1cm];
\draw plot [smooth, tension=1] coordinates {(0,-4) (0,-2)};
\draw plot [smooth, tension=1] coordinates {(0,-2) (3.46,-2)};
\draw plot [smooth, tension=0.9] coordinates {(0,-4) (-1,-3) (0,-1.4) (3.46,-2)};
\draw plot [smooth, tension=0.9] coordinates {(0,-4) (1.5,-3.5) (3.46,-2)};
\node[scale=0.7] at (0,-3){i};
\node[scale=0.7] at (1.4,-2){k'};
\node[scale=0.7] at (1.5,-3.5){1};
\node[scale=0.7] at (2.5,-1.75){j};
\end{tikzpicture} 
\end{equation}
where the arc labelled $1$ lives on the boundary. This follows by gluing certain vertices (\ref{diagonalflip}). In either case we remove the arc $k$ and replace it with the unique other arc $k'$ that forms a triangulation. For arcs inside self-folded triangles we do not allow mutation. Under these rules mutation of triangulations mimics exactly quiver mutation.
\begin{remark}
Due to the construction of quivers from triangulations, for any pair of arcs $i$ and $j$ forming a self-folded triangle (\ref{selffolded}) we have $b_{i,k}=b_{j,k}$ for all vertices $k$. So $\mu_i=\sigma\mu_j$ where $\sigma$ is the relabelling $i\leftrightarrow j$. Due to this, not being able to mutate at the arc inside the self-folded triangle, $i$, isn't really a restriction, we may as well just mutate at $j$. In \cite{clustersandtriangulated1} this inability to mutate inside self-folded triangles is rectified with the introduction of ``tagged arcs". 
\end{remark}
We also assign cluster variables to each arc of a triangulation. For this we let the boundary arcs have value $1$. To each arc $i$ in an initial triangulation we give initial cluster variables $x_i$. New cluster variables are defined to satisfy 
\begin{equation}\label{ptolemy}
x_{k}x_{k'}=x_ax_c+x_bx_d
\end{equation}
where the arcs are as labelled in (\ref{diagonalflip}). For quivers that are constructible as triangulated surfaces this matches the cluster mutation formula (\ref{clustermutationformula}). For the situation described in (\ref{selffoldedmutation}) we have
\[
x_{k}x_{k'}=1+x_j.
\]
Finally in the following situation
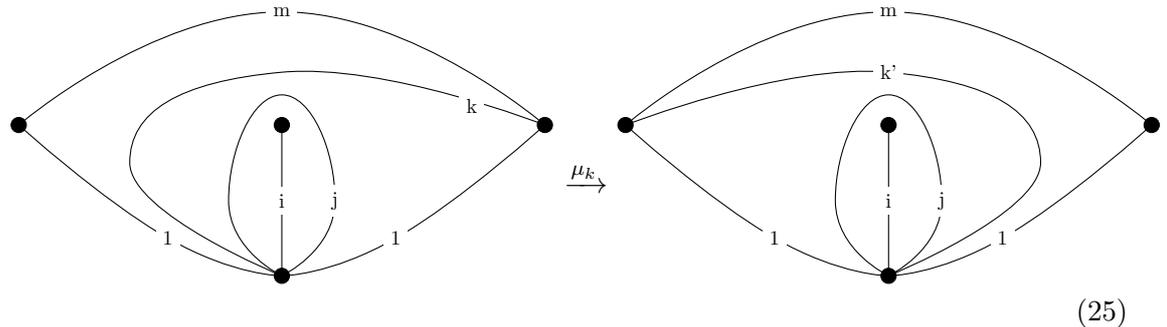
\begin{equation}
\begin{tikzpicture}[every node/.style={fill=white}]
\draw[fill=black] (0,-2) circle [radius=0.1cm];
\draw[fill=black] (0,-4) circle [radius=0.1cm];
\draw[fill=black] (3.46,-2) circle [radius=0.1cm];
\draw[fill=black] (-3.46,-2) circle [radius=0.1cm];
\draw plot [smooth, tension=1] coordinates {(0,-4) (0,-2)};
\draw plot [smooth, tension=1] coordinates {(0,-4) (-0.7,-3) (0,-1.6) (0.7,-3) (0,-4)};
\draw plot [smooth, tension=0.9] coordinates {(0,-4) (-1.5,-3.5) (-3.46,-2)};
\draw plot [smooth, tension=0.9] coordinates {(0,-4) (1.5,-3.5) (3.46,-2)};
\draw plot [smooth, tension=0.9] coordinates {(-3.46,-2) (0,-0.5) (3.46,-2)};
\draw plot [smooth, tension=0.9] coordinates {(0,-4) (-2,-2.5) (0,-1.3) (3.46,-2)};
\node[scale=0.7] at (0,-3){i};
\node[scale=0.7] at (0.7,-3){j};
\node at (4,-2.7){$\xrightarrow{\mu_k}$};
\node[scale=0.7] at (1.5,-3.5){1};
\node[scale=0.7] at (-1.5,-3.5){1};
\node[scale=0.7] at (0,-0.5){m};
\node[scale=0.7] at (2.5,-1.75){k};
\end{tikzpicture} 
\begin{tikzpicture}[every node/.style={fill=white}]
\draw[fill=black] (0,-2) circle [radius=0.1cm];
\draw[fill=black] (0,-4) circle [radius=0.1cm];
\draw[fill=black] (3.46,-2) circle [radius=0.1cm];
\draw[fill=black] (-3.46,-2) circle [radius=0.1cm];
\draw plot [smooth, tension=1] coordinates {(0,-4) (0,-2)};
\draw plot [smooth, tension=1] coordinates {(0,-4) (-0.7,-3) (0,-1.6) (0.7,-3) (0,-4)};
\draw plot [smooth, tension=0.9] coordinates {(0,-4) (-1.5,-3.5) (-3.46,-2)};
\draw plot [smooth, tension=0.9] coordinates {(0,-4) (1.5,-3.5) (3.46,-2)};
\draw plot [smooth, tension=0.9] coordinates {(-3.46,-2) (0,-0.5) (3.46,-2)};
\draw plot [smooth, tension=0.9] coordinates {(0,-4) (2,-2.5) (0,-1.3) (-3.46,-2)};
\node[scale=0.7] at (0,-3){i};
\node[scale=0.7] at (0.7,-3){j};
\node[scale=0.7] at (1.5,-3.5){1};
\node[scale=0.7] at (-1.5,-3.5){1};
\node[scale=0.7] at (0,-0.5){m};
\node[scale=0.7] at (0,-1.3){k'};
\end{tikzpicture} 
\end{equation}
we have $x_kx_{k'}=1+x_ix_jx_m$. Again, these are special cases of (\ref{ptolemy}) by gluing (\ref{diagonalflip}) in the right way.

We will use this construction of cluster algebras as triangulated surfaces to ``see" the cluster variables, in particular the frieze variables $X^i_n$ from (\ref{friezeformula}). We will also use this picture to show that the $A^m_n$ of (\ref{mainreccurence}) are indeed cluster variables. The recurrence defining these also means that they form a frieze, which we now define.
\subsection{Friezes}
In this section we give the definition of friezes, and then of repetition quivers. We show that in $\tilde{A}$ type the repetition quiver gives rise to a frieze. We define friezes on repetition quivers and show that for any acyclic quiver the frieze variables $X^i_n$ of (\ref{friezeformulaintro}) form a frieze on its repetition quiver.

Friezes were first defined in \cite{coxeterfrieze}. They were originally arrays of integers in the plane organised like brickwork and sandwiched between two rows of zeroes and ones:
\[\setcounter{MaxMatrixCols}{20}
\begin{matrix}
& \ldots & & 0 & & 0 & & 0 & & \ldots \\
\ldots & & 1 & & 1 & & 1 & & 1 & & \ldots \\
 & & \vdots & & \vdots & & \vdots & & \vdots & &  \\
& \ldots & & a & & b & & c & & \ldots \\
\ldots & & d & & e & & f & & g & & \ldots \\
& \ldots & & h & & i & & j & & \ldots \\
 & & \vdots & & \vdots & & \vdots & & \vdots & &  \\
\ldots & & 1 & & 1 & & 1 & & 1 & & \ldots \\
& \ldots & & 0 & & 0 & & 0 & & \ldots \\
\end{matrix}
\] 
such that each diamond 
$\begin{matrix}
& \beta & \\
\alpha & & \delta \\
& \gamma &
\end{matrix}$
satisfies $\alpha\delta-\beta\gamma=1$. For our purposes we don't demand that the entries are integers and don't require the rows of zeroes and ones. We have, for example, the friezes formed by the periodic quantities in $\tilde{A}$ and $\tilde{D}$ type shown in Proposition \ref{friezepropintro}, the proof of which we give in Section \ref{Friezeconstructionsection}. 

In order to see the link with the frieze patterns of Section \ref{dynamicalsystemssubsection}, we first need to make some definitions.
\begin{definition}\label{repetitionquiverdef}
For a quiver $Q$, the repetition $\mathbb{Z}Q$ is the quiver with vertices $(n,i)$, with $n\in \mathbb{Z}$ and $i$ a vertex of $Q$. For every arrow $i\mapsto j$ in $Q$ we have arrows $(n,i)\mapsto (n,j)$ and $(n-1,j)\mapsto (n,i)$ for every $n$. 
\end{definition}
The repetition quiver is an example of a stable translation quiver: it has a map $\tau:(n,i)\mapsto (n-1,i)$ such that there is a $1:1$ correspondence between arrows $j\mapsto i $ and $\tau(i)\mapsto j$.
\begin{definition}
A frieze on a repetition quiver is an assignment $f(n,i)$ for each vertex $(n,i)$ such that
\[
f(n,i)f(n-1,i)=1+\prod_{(m,j)\mapsto (n,i)} f(m,j)
\]
where the product is taken over all arrows into $(n,i)$.
\end{definition}
\begin{example}\label{repetitionexample}
We draw a frieze on the repetition quiver for $\tilde{A}_{5,2}$ (the quiver $\tilde{A}_{5,2}$ can be seen in red), where we identify the top and bottom rows:
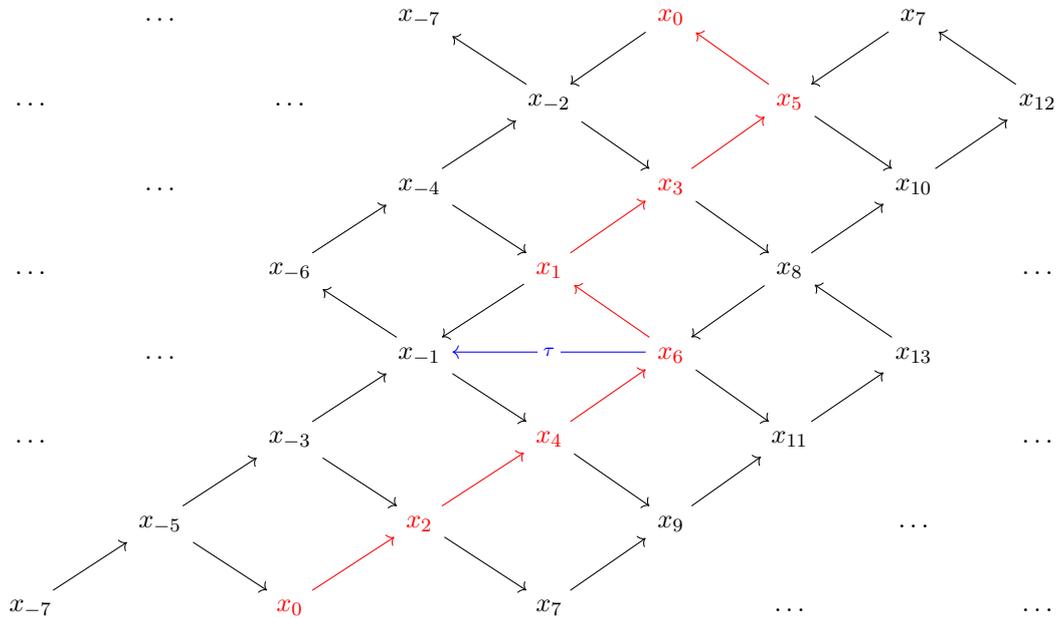
\[
\adjustbox{scale=0.9,center}{%
\begin{tikzcd}
&\ldots&&x_{-7}&& \color{red}x_0\arrow[dl] && x_7\arrow[dl]\\
\ldots&&\ldots&&x_{-2}\arrow[dr]\arrow[ul]&& \color{red}x_5\arrow[dr]\arrow[ul, red] &&x_{12}\arrow[ul]\\
&\ldots&&x_{-4}\arrow[dr]\arrow[ur]&&\color{red}x_3\arrow[dr] \arrow[ur, red] &&x_{10}\arrow[ur]\\
\ldots&&x_{-6}\arrow[ur]&& \color{red}x_1\arrow[dl] \arrow[ur, red] &&x_{8}\arrow[dl]\arrow[ur] && \ldots\\
&\ldots&&x_{-1}\arrow[dr]\arrow[ul]&& \color{red}x_6\arrow[ll, "\tau" description, blue]\arrow[dr]\arrow[ul, red] && x_{13}\arrow[ul]\\
\ldots&&x_{-3}\arrow[dr]\arrow[ur]&&\color{red}x_4 \arrow[dr]\arrow[ur, red] &&x_{11}\arrow[ur] &&\ldots\\
&x_{-5}\arrow[dr]\arrow[ur]&& \color{red}x_2\arrow[dr] \arrow[ur, red] &&x_9\arrow[ur] &&\ldots\\
x_{-7}\arrow[ur]&&\color{red}x_0 \arrow[ur, red] &&x_7\arrow[ur] && \ldots & & \ldots
\end{tikzcd}
}
\]
here the map $\tau$ acts by shifting each vertex left, as shown in blue for $\tau (x_6)=x_{-1}$. Furthermore this is a frieze due to the $\tilde{A}$ type relation (\ref{Atyperecurrence})
\[
x_{n}x_{n+7}=1+x_{n+2}x_{n+5}.
\]  
\end{example}
To define a frieze on the repetition quiver for $\tilde{A}_{q,p}$ for general $q$ and $p$ we let $f(n,i)=x_{nN+ip}$. The $x_n$ satisfy (\ref{Atyperecurrence}) so they give the frieze:
\begin{equation}\label{Atypefrieze}
\begin{matrix}
\ddots & & \vdots & & \vdots & & \iddots \\
&x_{n-q} & & x_{n+p} & & x_{n+N+p} & & \ldots \\
\ldots & & x_n & & x_{n+N} & & x_{n+2N} \\
&x_{n-p} && x_{n+q} & & x_{n+N+q} & & \ldots \\
\ldots &&x_{n+q-p} & & x_{n+2q} & & x_{n+N+2q} \\
& \iddots &&\vdots & & \vdots & & \ddots
\end{matrix}
\end{equation}
where, since 
\[
nN+ip=(n-p)N+(i+q+p)p
\] 
we identify $(n,i)=(n-p,i+q+p)$.
\begin{remark}
\begin{enumerate}[(i)]
\item The initial values for the associated dynamical system (\ref{Atyperecurrence}) are highlighted in red in Example \ref{repetitionexample}. In \cite{frises} this layout of the initial values is called a frontier (after rotation so that it is horizontal) and a formula is given for the frieze entries in terms of this frontier. 
\item  This type of construction can be found in \cite{calderochapoton} for $A$ type, with the initial cluster variables set to $1$, thus forming a traditional frieze.
\end{enumerate}
\end{remark}
For arbitrary quivers $Q$ (including $\tilde{D}$ type) we can take $f(n,i)=X^i_n$, as defined in (\ref{friezeformula}) to give a frieze on the repetition quiver $\mathbb{Z}Q$.
\subsection{Cluster algebras and representation theory}
Here we give the necessary background to compare our results with \cite{assemdupont}, where friezes are constructed on the various components of the cluster category. An overview of most of the elements of this section, including the cluster category and the cluster character, can be found in the notes \cite{kelleroverview}. Throughout this section the quiver $Q$ is taken to be acyclic and to have $N$ vertices and $\mathbf{k}$ is an arbitrary algebraically closed field.
\subsubsection{The cluster category}
In this section we define and briefly describe the Auslander-Reiten quivers for category of modules over the path algebra $\mathbf{k}Q$ if $Q$ is Dynkin or affine. We then describe the objects of the cluster category in these cases.
\begin{definition}
The Auslander-Reiten quiver $\Gamma(\mathbf{k}Q)$ of the module category $\mathrm{mod}\:\mathbf{k}Q$ has isomorphism classes of finitely generated, indecomposable modules for vertices. The arrows are given by ``irreducible maps". There is an automorphism 
\[
\tau:\Gamma(\mathbf{k}Q)\rightarrow \Gamma(\mathbf{k}Q)
\] 
called the Auslander-Reiten translate, such that $\tau(P)=0$ for all projectives $P$ and $\tau^{-1}(I)=0$ for all injectives $I$. All of the projective modules lie in the same component $\mathbf{p}$, which is called the preprojective component. Similarly the injectives all lie in $\mathbf{q}$, the preinjective component. The remaining components are called regular. See \cite{hugel} for an overview or \cite{auslander1974representation} for details.
\end{definition}
A quiver is called representation finite if, up to mutation, it has only finitely many indecomposable modules. Gabriel's theorem \cite{gabriel} says that these are precisely the Dynkin quivers. Moreover we have the following theorem, the first two points of which are from \cite{gabrielriedtmann} and paraphrased in \cite{hugel}. The third point can be found in \cite{crawleyboevey}.
\begin{theorem}
\begin{enumerate}[(i)]
\item If $Q$ is Dynkin then $\Gamma(\mathbf{k}Q)=\mathbf{p}=\mathbf{q}$ is a full and finite subquiver of $\mathbb{N}Q^{\mathrm{op}}$.
\item If $Q$ is not Dynkin then $\mathbf{p}=\mathbb{N}Q^{\mathrm{op}}$ and $\mathbf{q}=-\mathbb{N}Q^{\mathrm{op}}$ with $\mathbf{p}\cap \mathbf{q}=0$. The modules of $\mathbf{p}$ and $\mathbf{q}$ are uniquely determined by their dimension vectors. Finally $\mathbf{p}\cup \mathbf{q}\neq \Gamma(\mathbf{k}Q)$.
\item If $Q$ is affine then 
\[
\mathbf{p}=\{\tau^{-m}P_j\mid m\geq 0,\: j=1,\ldots, N\}, \qquad \mathbf{q}=\{\tau^{m}I_j\mid m\geq 0,\: j=1,\ldots, N\}
\]
and the regular modules $M$ are those that satisfy $\tau^m M\neq 0$ for all $m\in\mathbb{Z}$. The regular components are parametrised by $\lambda\in\mathbb{P}^1(\mathbf{k})$. By defining the quiver $A_{\infty}$ as 
\[
\begin{tikzpicture}[every node/.style={fill=white}]
\draw [->] (0,0) to (1.8,0);
\draw [->] (2,0) to (3.8,0);
\draw [->] (4,0) to (5.8,0);
\node at (6.2,0){\ldots};
\node at (0,0.0){$1$};
\node at (2,0.0){$2$};
\node at (4,0.0){$3$};
\end{tikzpicture}
\]
we can describe the regular components as ``tubes". They are given by ${\mathbb{Z}A_{\infty}/<\tau^{p_{\lambda}}>}$ for an appropriate $p_{\lambda}\in\mathbb{Z}_{\geq 0}$, known as the width of the tube, and ${\tau:(n,i)\mapsto(n-1,i)}$ where $n\in\mathbb{Z}$ and $i$ is a vertex of $A_{\infty}$.

The modules of the form $(n,1)$ for $n=1,2,\ldots,p_{\lambda}$ are called quasi-simple. If $p_{\lambda}=1$ then the tube is called homogeneous, otherwise it is called exceptional.  We define 
\[
\mathbb{P}^{\mathcal{E}}=\{\lambda\in\mathbb{P}^1(\mathbf{k})\mid p_{\lambda}>1\}
\]
and let $T_{\lambda}$ be the tube of width $p_{\lambda}$. A list of the widths of the exceptional tubes can be found in \cite{crawleyboevey}. In particular $|\mathbb{P}^{\mathcal{E}}|\leq 3$.
\end{enumerate}
\end{theorem}
The cluster category $\mathcal{C}_Q$ was first defined in \cite{tiltingcluster}. It has a suspension functor $[1]$ that coincides with the Auslander-Reiten translation $\tau$. The indecomposable objects of $\mathcal{C}_Q$ can be identified with the disjoint union of the objects of $\Gamma(\mathbf{k}Q)$ and the shifts of the projectives $P_i[1]$. The Auslander-Reiten quiver of $\mathcal{C}_Q$ is a stable translation quiver with translate $\tau$.

There is a connection between the cluster category and the cluster algebra given by the cluster character $X_?:\mathrm{Ob}(\mathcal{C}_Q)\rightarrow \mathcal{A}_Q$ which has a few components we need to define.
\subsubsection{The cluster character}
For a $\mathbf{k}Q$ module $M$ the quiver Grassmanian of dimension $\underline{e}\in\mathbb{Z}_{\geq 0}^N$ is 
\[
\mathrm{Gr}_{\underline{e}}(M)=\{N\subset M \mid \underline{\mathrm{dim}}(N)=\underline{e}\}
\]
and $\chi(\mathrm{Gr}_{\underline{e}}(M))$ is its Euler characteristic with respect to an appropriate cohomology. 

We let $b_{ij}$ be the number of arrows between vertices $i$ and $j$ in $Q$, then the Euler form $<-\:,\:->$ acts on pairs $a,a'\in \mathbb{Z}^N$ by
\[
<a\:,\:a'>=\sum_{i=1}^N a_ia'_i+\sum_{i,j=1}^N b_{ji}a_ia'_j.
\] 
\begin{definition}\cite{calderochapoton}
For a cluster algebra with initial cluster variables $u_i$, the cluster character is a map $X_?:\mathrm{Ob}(\mathcal{C}_Q)\rightarrow \mathcal{A}_Q$ defined by:
\begin{enumerate}[(i)]
\item If $M$ is an indecomposable $kQ$ module with $\underline{m}=\underline{\mathrm{dim}}(M)$ then
\[
X_M=\sum_{\underline{e}}\chi(\mathrm{Gr}_{\underline{e}}(M))\prod_i u_i^{-<\underline{e}\:,\:\underline{m}>-<\alpha_i\:,\:\underline{m}-\underline{e}>}
\]
where $\alpha_i=\underline{\mathrm{dim}}(S_i)$, the dimension of the simple module at $i$. The sum is taken over the $\underline{e}\in\mathbb{Z}^N$ such that $\chi(\mathrm{Gr}_{\underline{e}}(M))\neq 0$ and the product is taken over all vertices $i$. 
\item If $M=P_i[1]$ then $X_M=u_i$.
\item For any two objects $M$ and $N$
\[
X_{M\bigoplus N}=X_MX_N.
\] 
\end{enumerate}
\end{definition}
\begin{theorem} \cite{calderokellerii} Theorem 4. The map $X_?$ gives a bijection between the set of isomorphism classes of rigid indecomposable modules in $\mathcal{C}_Q$ and the set of cluster variables $\mathcal{A}_Q$.  
\end{theorem}
\begin{proposition}\cite{assemdupont} Proposition 2.2. For an acyclic quiver $Q$ the cluster character $X_?$ induces a frieze on the repetition quiver $\Gamma(\mathcal{C}_Q)$.
\end{proposition}
\subsubsection{Friezes on the cluster category}
We are interested affine quivers, in which case the transjective component of $\Gamma(\mathcal{C}_Q)$ is isomorphic to the repetition quiver $\mathbb{Z}Q$ and has the objects $P_i[1]$ at the vertices $(0,i)$. The frieze on $\mathbb{Z}Q$ given by (\ref{friezeformula}) places the initial cluster variables $X^i_0$ at $(0,i)$, so by \cite{assemdupont}, Corollary 3.2, the frieze $X^i_n$ coincides with the frieze given by $X_?$ on the transjective component.

Throughout this section the vertex $e$ is taken to be a fixed extending vertex. In the $\tilde{A}$ case the extending vertices are all vertices and in the $\tilde{D}$ case these are the vertices labelled $1,2,N$ or $N+1$ in Figure \ref{Dquiver}. When we use these results in Sections \ref{Atypeclusterfriezesection} and \ref{Dtypeclusterfriezebit} we choose $e$ explicitly.

For any $\lambda\in\mathbb{P}^1(k)$ there exists a unique quasi-simple module $M_{\lambda}$ in $\mathcal{T}_{\lambda}$ such that $\mathrm{dim}\:M_{\lambda}(e)=1$. Set $N_{\lambda}=M_{\lambda}[1]$ if $e$ is a source or $N_{\lambda}=M_{\lambda}[-1]$ if $e$ is a sink. In \cite{assemdupont} it is proved that there exists transjective $B_{\lambda}$ and  $B'_{\lambda}$ that are uniquely determined by the existence of non-split triangles 
\[
\begin{tikzcd}
N_{\lambda}\arrow[r] & B_{\lambda}\arrow[r] & S_{e}\arrow[r] & N_{\lambda}[1]
\end{tikzcd}
\qquad
\begin{tikzcd}
S_e\arrow[r] & B'_{\lambda}\arrow[r] & N_{\lambda}\arrow[r] & S_e[1]
\end{tikzcd}
\]
for a source $e$ and 
\[
\begin{tikzcd}
N_{\lambda}\arrow[r] & B'_{\lambda}\arrow[r] & S_{e}\arrow[r] & N_{\lambda}[1]
\end{tikzcd}
\qquad
\begin{tikzcd}
S_e\arrow[r] & B_{\lambda}\arrow[r] & N_{\lambda}\arrow[r] & S_e[1]
\end{tikzcd}
\]
if $e$ is a sink. For $\tilde{A}$ and $\tilde{D}$ type these $B_{\lambda}$ and $B'_{\lambda}$ are then constructed explicitly. These triangles are used to allow consideration of friezes on the whole of $\Gamma(\mathcal{C}_Q)$.
\begin{definition}
For an affine quiver $Q$, a cluster frieze on $\Gamma(\mathcal{C}_{Q})$ is a frieze $f$ on $\Gamma(\mathcal{C}_{Q})$ such that, for any $\lambda\in\mathbb{P}^{\mathcal{E}}$ we have
\begin{equation}\label{clusterfrieze}
f(N_{\lambda}[k])=\frac{f(B_{\lambda}[k])+f(B'_{\lambda}[k])}{f(S_{e}[k])}
\end{equation}
for any $k=1,2,\ldots p_{\lambda}-1$.
\end{definition}
We can consider this as ``gluing" the friezes on the exceptional tubes to the frieze on the transjective component. Due to (\ref{clusterfrieze}) the values of a cluster frieze on the exceptional tubes are determined by the values on the transjective component. 
\begin{proposition}\cite{assemdupont}, Proposition 2.2. If $f$ is a cluster frieze on $\Gamma(\mathcal{C}_Q)$ such that $f(P_i[1])=u_i$ for each $i$ then $f(M)=X_M$ for all $M\in\mathcal{C}_Q$.
\end{proposition}
For $\tilde{A}$ and $\tilde{D}$ quivers we have seen that the $X^i_j$ of (\ref{friezeformula}) form a frieze on their transjective components.  Our goal is to define friezes on the exceptional components, in terms of the periodic quantities $J_n$ and $\tilde{J}_n$ ($\tilde{A}$ type) and $J'_n$ ($\tilde{D}$ type) satisfying the cluster frieze condition (\ref{clusterfrieze}). This is done in Sections \ref{Atypeclusterfriezesection} and \ref{Dtypeclusterfriezebit}.
\section{Linear relations and friezes for the matrix entries $A^m_n$, $\tilde{A}^m_n$ and ${A'}^m_n$}\label{Friezeconstructionsection}
In this section we prove Propositions \ref{matrixentryclustervariables} and \ref{Dtypematrixentries} which give linear relations between the $A^m_n$ and $\tilde{A}^m_n$ in the $\tilde{A}$ case and ${A'}^m_n$ in the $\tilde{D}$ case, the matrix entries that appear in the construction of the constant coefficient linear relations. For both we show that we can construct friezes from these matrix entries.

We describe the entries of $M_n^m$ as
\[
M^m_n:=
\begin{pmatrix}
A^m_n & B^m_n \\
C^m_n & D^m_n
\end{pmatrix}
\]
These satisfy $M^{m+1}_n=M^m_nL_{n+mp}$ so we have
\[
\begin{pmatrix}
A^{m+1}_{n} & B^{m+1}_{n} \\
C^{m+1}_{n} & D^{m+1}_{n}
\end{pmatrix}
=
\begin{pmatrix}
A^m_n & B^m_n \\
C^m_n & D^m_n
\end{pmatrix}
\begin{pmatrix}
J_{n+mp} & 1\\
-1 & 0
\end{pmatrix}
=
\begin{pmatrix}
A^m_nJ_{n+mp}-B^m_n & A^m_n \\
C^m_nJ_{n+mp}-D_n & C^m_n
\end{pmatrix}
\]
This gives $B^m_{n}=A^{m-1}_{n}$ and $D^m_n=C^{m-1}_{n}$, so $A^m_n$ and $C^m_n$ satisfy
\[
A^{m+1}_{n}=A^m_nJ_{n+mp}-A^{m-1}_{n}, \qquad C^{m+1}_{n}=C^m_nJ_{n+mp}-C^{m-1}_{n}
\]
with initial values
\[
A^1_n=J_n, \qquad A^0_n=1, \qquad C^1_n=-1, \qquad C^0_n=0
\] 
since $M^0_n=I$ and $M^1_n=L_n$. We have $C^2_n=-J_{n+p}=-A^1_{n+p}$ and $C^1_n=-1=-A^0_{n+p}$ so $C^m_n=-A^{m-1}_{n+p}$ for all $n,m$. This proves Proposition \ref{matrixentryclustervariables} for $M^m_n$ and the result for $\tilde{M}^m_n$ is proved similarly.
\begin{remark}
In \cite{fordyhone}, for $q=1$, the traces $\mathcal{K}^m_n=A^m_n+D^m_n$ are constructed in terms of a recursion operator 
\[
\mathcal{R}^{(m)}_n=J_{n+m}J_{n+m+1}\frac{\partial^2}{\partial J_n\partial J_{n+m-1}}-J_{n+m}\frac{\partial}{\partial J_n}-J_{n+m+1}\frac{\partial}{\partial J_{n+m-1}}+J_{n+m}J_{n+m+1}-1
\]
satisfying $\mathcal{K}^{m+2}_n=\mathcal{R}^m_n\mathcal{K}^m_n$.
\end{remark}
The recurrence (\ref{mainreccurence}) appears in \cite{coxeterfrieze}, with the iterates forming a frieze with the boundaries of ones and zeroes. Coxeter gives a general solution, which still works in our case; we define the determinants
\begin{equation}\label{determinant}
{D}^{m}_{a}(F_n):=
\begin{vmatrix}
F_n & 1 & 0 \\
1 & F_{n+a} & 1 & 0 \\
0 & 1 & F_{n+2a} & 1 & 0 \\
& 0 & 1 & F_{n+3a} & \ddots \\
& & 0 & \ddots & \ddots & 1 \\
& & & & 1 & F_{n+(m-1)a}
\end{vmatrix}
\end{equation}
for an arbitrary function $F_n$ and integer $a$. By expanding along the last row these satisfy
\[
{D}^{m}_{a}(F_n)=F_{n+(m-1)a}{D}^{m-1}_{a}(F_n)-{D}^{m-2}_{a}(F_n)
\]
with ${D}^{0}_{a}(F_n)=1$ and ${D}^{1}_{a}(F_n)=F_n$, hence we can express each of $A^m_n$, $\tilde{A}^m_n$ and ${A'}^m_n$ in terms of these determinants:
\[
{D}^{m}_{p}(J_n)=A^m_n, \qquad D^m_q(\tilde{J}_n)=\tilde{A}^m_n, \qquad D^m_1(J'_n)={A'}^m_n.
\]
Furthermore applying the Desnanot-Jacobi identity to ${D}^{m}_{a}(F_n)$ gives 
\[
D^{m}_a(F_n)D^{m-2}_a(F_{n+a})=D^{m-1}_a(F_n)D^{m-1}_a(F_{n+a})-1
\] 
so the determinants $D^m_a(F_n)$ form a frieze:
\begin{equation}\label{generaldetfrieze}\setcounter{MaxMatrixCols}{20}
\begin{matrix}
&1 & & 1 & & 1 & & 1 & & \ldots \\
\ldots & & F_n & & F_{n+a} & & F_{n+2a} & & F_{n+3a} \\
&D^2_a(F_{n-a}) & & D^2_a(F_n) & & D^2_a(F_{n+a}) & & D^2_a(F_{n+2a}) & & \ldots \\
\ldots && D^3_a(F_{n-a}) & & D^3_a(F_{n}) & & D^3_a(F_{n+a}) & & D^3_a(F_{n+2a}) \\
&\vdots & & \vdots & & \vdots && \vdots 
\end{matrix}
\end{equation}
This gives two friezes in $\tilde{A}$ type:
\begin{equation}\label{frieze}
(i) \qquad \:\: F_n=J_n, \qquad a=p.
\end{equation}
\begin{equation}\label{tildefrieze}
(ii) \qquad F_n=\tilde{J}_n, \qquad a=q.
\end{equation}
and one frieze in $\tilde{D}$ type, 
\begin{equation}\label{friezeDtype}
\qquad \quad \:\: F_n=J'_n, \qquad a=1.
\end{equation}
\begin{remark}
\begin{enumerate}[(i)]
\item The determinants (\ref{determinant}) are examples of continuants. See \cite{muir}, for example, where they are shown to be related to continued fractions.
\item This determinant solution for Dynkin $A$ type appears in \cite{dupont}. There all cluster variables $x_n$ come from the cluster map and they form a frieze (\ref{frieze}), with the $x_n$ replacing the $J_n$ in an appropriate way. There the $D^m_n$ are called ``generalised Chebyshev polynomials".
\end{enumerate}
\end{remark}
\section{Triangulated surfaces and the cluster frieze in $\tilde{A}$ type}\label{Atypesection}
In Section \ref{Atypetriangulatedsection} we construct $\tilde{A}_{q,p}$ cluster algebras as triangulated surfaces. We first find the cluster variables $x_n$, obtained by the cluster map, as arcs before proving that the remaining arcs are given by the determinants $D^l_p(J_{jp})$ and $D^l_q(\tilde{J}_{jq})$ of Theorem \ref{findingallclustervarstheorem}, proving the $\tilde{A}$ part of that theorem. Finally in Section \ref{Atypeclusterfriezesection} we prove that the friezes constructed in this case, (\ref{Atypefrieze}) for the $x_n$ and (\ref{frieze}) and (\ref{tildefrieze}) for the $J_n$ and $\tilde{J}_n$, are precisely the cluster friezes given in \cite{assemdupont}. 
\subsection{$\tilde{A}$ type cluster algebras as triangulated surfaces}\label{Atypetriangulatedsection}
The $\tilde{A}_{q,p}$ cluster algebra is obtained by triangulations of an annulus with $q$ and $p$ vertices on the inner and outer boundaries, respectively. By cutting this annulus we represent the surface as a  strip with $q$ vertices on the top and $p$ vertices on the bottom, as in Figure \ref{annuluswithsomearcs}, where we have displayed some possible arcs (cluster variables) but not a full triangulation.  

\begin{figure}
\centering
\begin{tikzpicture}
\draw[fill=black] (0,0) circle [radius=0.1cm];
\draw[fill=black] (2,0) circle [radius=0.1cm];
\draw[fill=black] (4,0) circle [radius=0.1cm];
\draw[fill=black] (3,2) circle [radius=0.1cm];
\draw[fill=black] (5,2) circle [radius=0.1cm];
\draw[fill=black] (7,2) circle [radius=0.1cm];
\draw plot [smooth, tension=1] coordinates {(0,0) (2,0.7) (4,0)};
\draw plot [smooth, tension=1] coordinates {(3,2) (5,1.3) (7,2)};
\draw plot [smooth, tension=1] coordinates {(0,0) (2,0)};
\draw plot [smooth, tension=1] coordinates {(0,2) (7,2)};
\draw plot [smooth, tension=1] coordinates {(2,0) (7,0)};
\node at (2,1.0) {$J_{jp}$};
\node at (5,1.0) {$\tilde{J}_{-iq}$};
\node at (0,-0.3) {$j$};
\node at (2,-0.3) {$j+1$};
\node at (4,-0.3) {$j+2$};
\node at (3,2.3) {$q+i-2$};
\node at (5,2.3) {$q+i-1$};
\node at (7,2.3) {$q+i$};
\end{tikzpicture}\caption{The arcs $J_{jp}$ and $\tilde{J}_{-iq}$.}\label{arcJ}
\end{figure}
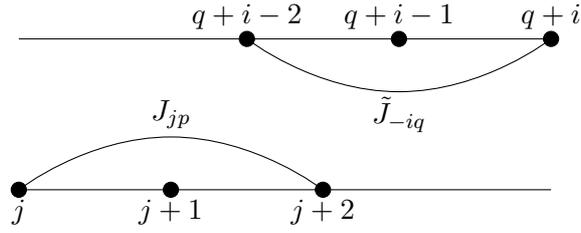

In order to present the $\tilde{A}_{q,p}$ quiver constructed in Section \ref{dynamicalsystemssubsection} (for example Figure \ref{EuclideanA15}) as a triangulation, we take a strip with vertices labelled as in Figure  \ref{annuluswithsomearcs} and draw the arcs $kp$ for ${k=0,1,\ldots,p+q-2}$ from vertex $k- \lfloor{kp/N} \rfloor$ (on the bottom) to vertex $q+\left \lfloor{kp/N}\right \rfloor$ (on the top). We also have an arc labelled $(p+q-1)p$ from $0$ on the bottom to $p+q-1$ on the top. We then reduce the labels modulo $N$.
\begin{example}
In Figure \ref{exampletriangulation} we demonstrate this construction for the quiver from Figure \ref{EuclideanA15}. Before reducing modulo $N$ our vertex labels are $0,7,14,21,28,\ldots, 98$. The first three, $k=0,1,2$, satisfy $\lfloor{kp/N} \rfloor=0$ so they give arcs $k\mapsto q=8$ labelled $0,7$ and $14$. The next, $k=3$ has $\lfloor{kp/N} \rfloor=1$ so it gives an arc $2\mapsto q+1=9$ which we label $6$, after reducing $21$ modulo $15$. Continuing this gives the rest of the triangulation.
\begin{figure}
\centering
\[
\begin{tikzpicture}[every node/.style={fill=white}]
\foreach \x in {0,...,7}
\draw[fill=black] (\x,0) circle [radius=0.1cm];
\foreach \x in {0,...,6}
\draw[fill=black] (\x+1/2,2) circle [radius=0.1cm];
\foreach \x in {-2,...,8}
\draw [-] (\x,0) to (\x+1,0);
\foreach \x in {-2,...,8}
\draw [-] (\x,2) to (\x+1,2);
\foreach \x in {0,...,7}
\node[scale=0.7] at (\x,-0.3) {\x};
\node[scale=0.7] at (0.5,2.3) {$8$};
\node[scale=0.7] at (1.5,2.3) {$9$};
\node[scale=0.7] at (2.5,2.3) {$10$};
\node[scale=0.7] at (3.5,2.3) {$11$};
\node[scale=0.7] at (4.5,2.3) {$12$};
\node[scale=0.7] at (5.5,2.3) {$13$};
\node[scale=0.7] at (6.5,2.3) {$14$};
\node[scale=0.7] at (8,-0.3) {$0$};
\node[scale=0.7] at (-1.5,2.3) {$14$};
\draw [-] (0,0) to node[scale=0.7] {0} (0.5,2);
\draw [-] (1,0) to node[scale=0.7] {7} (0.5,2);
\draw [-] (2,0) to node[scale=0.7] {14} (0.5,2);
\draw [-] (2,0) to node[scale=0.7] {6} (1.5,2);
\draw [-] (3,0) to node[scale=0.7] {13} (1.5,2);
\draw [-] (3,0) to node[scale=0.7] {5} (2.5,2);
\draw [-] (4,0) to node[scale=0.7] {12} (2.5,2);
\draw [-] (4,0) to node[scale=0.7] {4} (3.5,2);
\draw [-] (5,0) to node[scale=0.7] {11} (3.5,2);
\draw [-] (5,0) to node[scale=0.7] {3} (4.5,2);
\draw [-] (6,0) to node[scale=0.7] {10} (4.5,2);
\draw [-] (6,0) to node[scale=0.7] {2} (5.5,2);
\draw [-] (7,0) to node[scale=0.7] {9} (5.5,2);
\draw [-] (7,0) to node[scale=0.7] {1} (6.5,2);
\draw [-] (8,0) to node[scale=0.7] {8} (6.5,2);
\draw [-] (0,0) to node[scale=0.7] {8} (-1.5,2);
\draw[dotted] (-0.5,-1) to (-0.5,3);
\draw[dotted] (7.5,-1) to (7.5,3);
\draw[fill=black] (-1.5,2) circle [radius=0.1cm];
\draw[fill=black] (8,0) circle [radius=0.1cm];
\end{tikzpicture}
\]
\caption{The triangulation for the $\tilde{A}_{8,7}$ quiver.}
\label{exampletriangulation}
\end{figure}
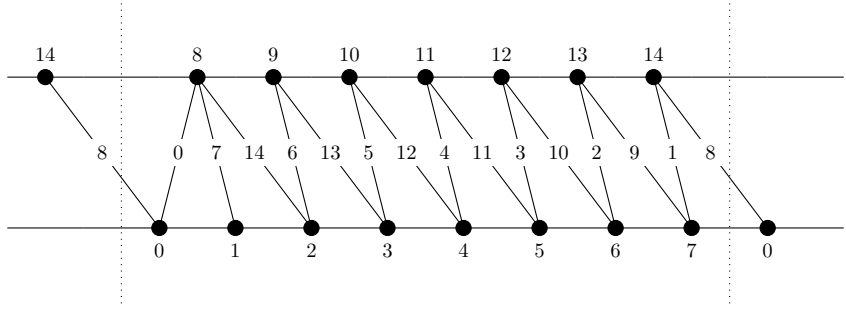
\end{example}
We see that in this construction every arc connects the top boundary to the bottom. Of course, any arc will either do this or connect two vertices on the same boundary component. The following proposition tells us precisely how to obtain every arc of the first kind.
\begin{proposition}\label{Atypebipartitebelt}
The arcs connecting the top of the strip to the bottom are in bijection with the cluster variables $x_n$ obtained by the recurrence (\ref{Atyperecurrence}).
\begin{proof}
We temporarily forget the periodicity of the strip and consider the (now infinitely many) vertices on the top and bottom to be labelled by $v_i$ and $v'_i$ for $i\in\mathbb{Z}$, as in Figure \ref{alternativelabelling}.
\begin{figure}
\centering
\begin{tikzpicture}
\draw[fill=black] (0,0) circle [radius=0.1cm];
\draw[fill=black] (1,0) circle [radius=0.1cm];
\draw[fill=black] (2,0) circle [radius=0.1cm];
\node at (3,0) {$\ldots$};
\draw[fill=black] (4,0) circle [radius=0.1cm];
\draw[fill=black] (0.5,2) circle [radius=0.1cm];
\draw[fill=black] (1.5,2) circle [radius=0.1cm];
\node at (2.5,2) {$\ldots$};
\draw[fill=black] (3.5,2) circle [radius=0.1cm];
\foreach \x in {-1,...,1}
\draw [-] (\x,0) to (\x+1,0);
\draw [-] (2,0) to (2.5,0);
\draw [-] (3.5,0) to (4,0);
\foreach \x in {4,...,6}
\draw [-] (\x,0) to (\x+1,0);
\foreach \x in {-1,...,1}
\draw [-] (\x,2) to (\x+1,2);
\foreach \x in {3,...,6}
\draw [-] (\x,2) to (\x+1,2);
\draw [-] (-1.5,2) to (-0.5,2);
\node at (0,-0.3) {$v_0$};
\node at (1,-0.3) {$v_1$};
\node at (2,-0.3) {$v_2$};
\node at (4,-0.3) {$v_{q-1}$};
\node at (5,-0.3) {$v_{q}$};
\node at (6,-0.3) {$v_{q+1}$};
\node at (7,-0.3) {$v_{q+2}$};
\node at (-1.5,2.3) {$v'_{-1}$};
\node at (0.5,2.3) {$v'_0$};
\node at (1.5,2.3) {$v'_{1}$};
\node at (3.5,2.3) {$v'_{p-1}$};
\node at (5.5,2.3) {$v'_{p}$};
\node at (6.5,2.3) {$v'_{p+1}$};
\draw[fill=black] (5,0) circle [radius=0.1cm];
\draw[fill=black] (6,0) circle [radius=0.1cm];
\draw[fill=black] (7,0) circle [radius=0.1cm];
\draw[fill=black] (5.5,2) circle [radius=0.1cm];
\draw[fill=black] (-1.5,2) circle [radius=0.1cm];
\draw[fill=black] (6.5,2) circle [radius=0.1cm];
\draw[dotted] (-0.5,-1) to (-0.5,3);
\draw[dotted] (4.5,-1) to (4.5,3);
\end{tikzpicture}\caption{Alternative labelling of the vertices on the strip.}\label{alternativelabelling}
\end{figure}
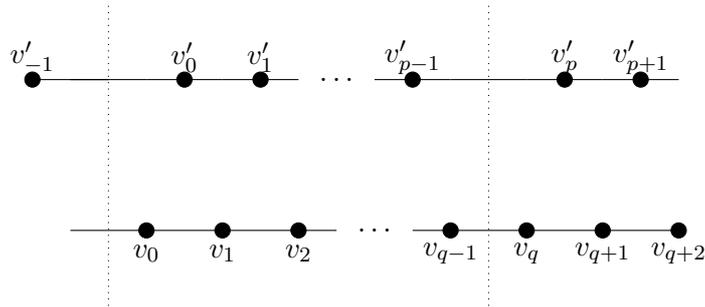
The mutation $\mu_{k}$ will mutate each arc $k$ while it is a diagonal of a square (\ref{DiagonalAtype})
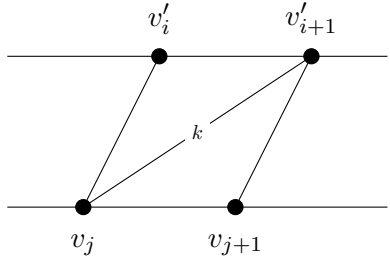
\begin{equation}\label{DiagonalAtype}
\begin{tikzpicture}[every node/.style={fill=white}]
\draw[fill=black] (0,0) circle [radius=0.1cm];
\draw[fill=black] (2,0) circle [radius=0.1cm];
\draw[fill=black] (1,2) circle [radius=0.1cm];
\draw[fill=black] (3,2) circle [radius=0.1cm];
\draw [-] (-1,0) to (4,0);
\draw [-] (-1,2) to (4,2);
\draw [-] (0,0) to node[scale=0.7] {$k$} (3,2);
\draw [-] (2,0) to (3,2);
\draw [-] (0,0) to (1,2);
\node[scale=1] at (0,-0.5) {$v_j$};
\node[scale=1] at (2,-0.5) {$v_{j+1}$};
\node[scale=1] at (1,2.5) {$v'_i$};
\node[scale=1] at (3,2.5) {$v'_{i+1}$};
\end{tikzpicture}
\end{equation}
to give the arc on the other diagonal $v'_{i}\:\mbox{---}\:v_{j+1}$, for example mutating the triangulation of Figure \ref{exampletriangulation} at $0$ gives Figure \ref{triangulationaftermutation}. The whole mutation sequence ${\mu=\mu_N\mu_{N-1}\ldots\mu_1\mu_0}$ will do this to every arc.
 
With the labelling in Figure \ref{alternativelabelling}, the initial quiver has arcs $kp$ for ${k=0,1,\ldots,p+q-1}$ from vertex $k- \lfloor{kp/N} \rfloor$ (on the bottom) to vertex $\left \lfloor{kp/N}\right \rfloor$ (on the top).

To show that a general arc $v_{i}\:\mbox{---}\:v'_{j}$ appears due to the cluster map, we reinstate the periodicity on Figure \ref{alternativelabelling} and have that 
\[
v_{i}\:\mbox{---}\:v'_{j}=v_{i+mp}\:\mbox{---}\:v'_{j+mq}
\]
for all $m\in\mathbb{Z}$, so by taking the correct $m$ we can assume that $0\leq i+j\leq p+q-1$. In this case the arc $k=i+j$ that appears in the initial triangulation is from $k- \lfloor{kp/N} \rfloor$ on the bottom to $\lfloor{kp/N}\rfloor $ on the top. Applying $\mu^{\lfloor{kp/N}\rfloor-j}$ will give the arc $v_{i}\:\mbox{---}\:v'_{j}$.
\end{proof}
\end{proposition}
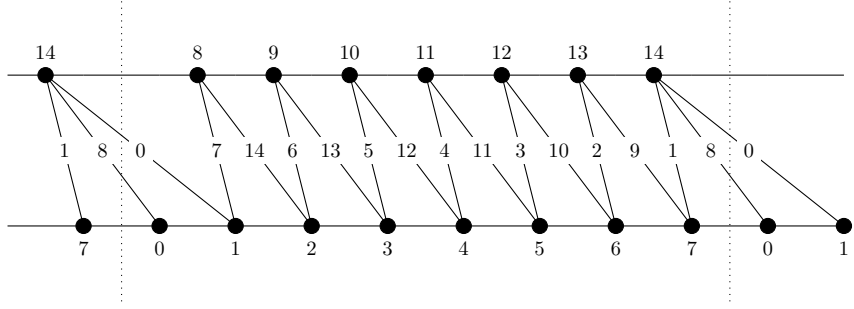
\begin{figure}
\centering
\[
\begin{tikzpicture}[every node/.style={fill=white}]
\foreach \x in {0,...,7}
\draw[fill=black] (\x,0) circle [radius=0.1cm];
\foreach \x in {0,...,6}
\draw[fill=black] (\x+1/2,2) circle [radius=0.1cm];
\foreach \x in {-2,...,8}
\draw [-] (\x,0) to (\x+1,0);
\foreach \x in {-2,...,8}
\draw [-] (\x,2) to (\x+1,2);
\foreach \x in {0,...,7}
\node[scale=0.7] at (\x,-0.3) {\x};
\node[scale=0.7] at (0.5,2.3) {$8$};
\node[scale=0.7] at (1.5,2.3) {$9$};
\node[scale=0.7] at (2.5,2.3) {$10$};
\node[scale=0.7] at (3.5,2.3) {$11$};
\node[scale=0.7] at (4.5,2.3) {$12$};
\node[scale=0.7] at (5.5,2.3) {$13$};
\node[scale=0.7] at (6.5,2.3) {$14$};
\node[scale=0.7] at (8,-0.3) {$0$};
\node[scale=0.7] at (9,-0.3) {$1$};
\node[scale=0.7] at (-1,-0.3) {$7$};
\node[scale=0.7] at (-1.5,2.3) {$14$};
\draw [-] (1,0) to node[scale=0.7] {0} (-1.5,2);
\draw [-] (1,0) to node[scale=0.7] {7} (0.5,2);
\draw [-] (2,0) to node[scale=0.7] {14} (0.5,2);
\draw [-] (2,0) to node[scale=0.7] {6} (1.5,2);
\draw [-] (3,0) to node[scale=0.7] {13} (1.5,2);
\draw [-] (3,0) to node[scale=0.7] {5} (2.5,2);
\draw [-] (4,0) to node[scale=0.7] {12} (2.5,2);
\draw [-] (4,0) to node[scale=0.7] {4} (3.5,2);
\draw [-] (5,0) to node[scale=0.7] {11} (3.5,2);
\draw [-] (5,0) to node[scale=0.7] {3} (4.5,2);
\draw [-] (6,0) to node[scale=0.7] {10} (4.5,2);
\draw [-] (6,0) to node[scale=0.7] {2} (5.5,2);
\draw [-] (7,0) to node[scale=0.7] {9} (5.5,2);
\draw [-] (7,0) to node[scale=0.7] {1} (6.5,2);
\draw [-] (8,0) to node[scale=0.7] {8} (6.5,2);
\draw [-] (0,0) to node[scale=0.7] {8} (-1.5,2);
\draw [-] (-1,0) to node[scale=0.7] {1} (-1.5,2);
\draw [-] (6.5,2) to node[scale=0.7] {0} (9,0);
\draw[dotted] (-0.5,-1) to (-0.5,3);
\draw[dotted] (7.5,-1) to (7.5,3);
\draw[fill=black] (-1.5,2) circle [radius=0.1cm];
\draw[fill=black] (8,0) circle [radius=0.1cm];
\draw[fill=black] (9,0) circle [radius=0.1cm];
\draw[fill=black] (-1,0) circle [radius=0.1cm];
\end{tikzpicture}
\]
\caption{The triangulation for the $\tilde{A}_{8,7}$ quiver after mutation at $0$.}
\label{triangulationaftermutation}
\end{figure}
Aside from the arcs discussed in the previous proposition, we only have arcs connecting the top (or bottom) of the strip to itself, for example the arc $2\:\mbox{---}\: p$ in Figure \ref{annuluswithsomearcs}. We call these $J_{jp}$ and $\tilde{J}_{-iq}$, as shown in Figure \ref{arcJ}, for $j=0,1,\ldots,q-1$ and $i=0,1,\ldots,p-1$. 
This naming is justified by the following lemma.
\begin{lemma}\label{constructingJarcslemma}
The arcs $J_{jp}$ and $\tilde{J}_{-iq}$ are the arcs associated with the periodic quantities (\ref{periodicquantities}), which are cluster variables. 
\begin{proof}
We first prove that for any arc labelled $\overline{kp}$, for some $k$, at vertex $j$ on the bottom boundary we have $\overline{kp}\equiv jp \mod q$. Here $\overline{kp}$ denotes the label after reduction modulo $N$. We prove this by induction and assume that the statement is true for the arcs labelled up to $\overline{(k-1)p}$. There are two possibilities:
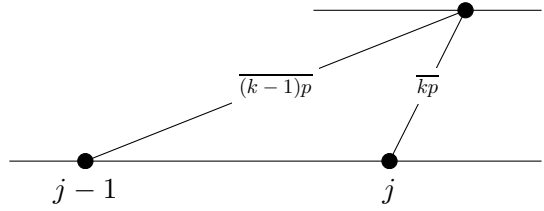
\begin{equation*}
\begin{tikzpicture}[every node/.style={fill=white}]
\draw[fill=black] (0,0) circle [radius=0.1cm];
\draw[fill=black] (4,0) circle [radius=0.1cm];
\draw[fill=black] (5,2) circle [radius=0.1cm];
\draw [-] (-1,0) to (6,0);
\draw [-] (3,2) to (6,2);
\draw [-] (4,0) to node[scale=0.7] {$\overline{kp}$} (5,2);
\draw [-] (0,0) to node[scale=0.7] {$\overline{(k-1)p}$} (5,2);
\node at (0,-0.4) {$j-1$};
\node at (4,-0.4) {$j$};
\end{tikzpicture}
\end{equation*}
or
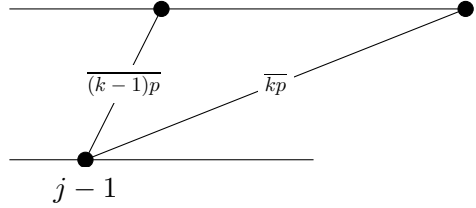
\begin{equation*}
\begin{tikzpicture}[every node/.style={fill=white}]
\draw[fill=black] (4,0) circle [radius=0.1cm];
\draw[fill=black] (9,2) circle [radius=0.1cm];
\draw[fill=black] (5,2) circle [radius=0.1cm];
\draw [-] (3,0) to (7,0);
\draw [-] (3,2) to (9,2);
\draw [-] (4,0) to node[scale=0.7] {$\overline{(k-1)p}$} (5,2);
\draw [-] (4,0) to node[scale=0.7] {$\overline{kp}$} (9,2);
\node at (4,-0.4) {$j-1$};
\end{tikzpicture}
\end{equation*}
In the first case $\overline{kp}=\overline{(k-1)p}+p$ so 
\[
\overline{kp}\equiv \overline{(k-1)p}+p\equiv (j-1)p+p\equiv jp \mod q
\]
where we have used the induction assumption. In the second case 
\[
\overline{kp}=\overline{(k-1)p}+p-N=\overline{(k-1)p}-q
\] 
so 
\[
\overline{kp}\equiv \overline{(k-1)p} \equiv (j-1)p \mod q.
\] 
To construct the $J_n$ we let $\overline{kp}$ be the rightmost arc at $j$. Near $j+1$ there are two situations:
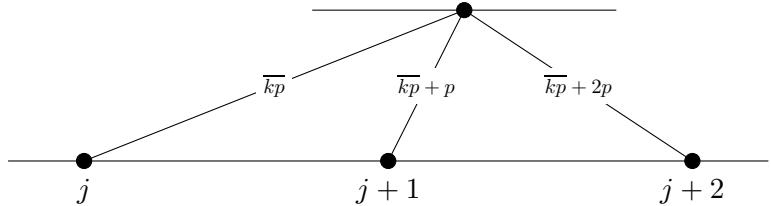
\begin{equation*}
\begin{tikzpicture}[every node/.style={fill=white}]
\draw[fill=black] (0,0) circle [radius=0.1cm];
\draw[fill=black] (4,0) circle [radius=0.1cm];
\draw[fill=black] (8,0) circle [radius=0.1cm];
\draw[fill=black] (5,2) circle [radius=0.1cm];
\draw [-] (-1,0) to (9,0);
\draw [-] (3,2) to (7,2);
\draw [-] (4,0) to node[scale=0.7] {$\overline{kp}+p$} (5,2);
\draw [-] (0,0) to node[scale=0.7] {$\overline{kp}$} (5,2);
\draw [-] (8,0) to node[scale=0.7] {$\overline{kp}+2p$} (5,2);
\node at (0,-0.4) {$j$};
\node at (4,-0.4) {$j+1$};
\node at (8,-0.4) {$j+2$};
\end{tikzpicture}
\end{equation*}
or
\begin{equation*}
\begin{tikzpicture}[every node/.style={fill=white}]
\draw[fill=black] (0,0) circle [radius=0.1cm];
\draw[fill=black] (4,0) circle [radius=0.1cm];
\draw[fill=black] (8,0) circle [radius=0.1cm];
\draw[fill=black] (9,2) circle [radius=0.1cm];
\draw[fill=black] (5,2) circle [radius=0.1cm];
\draw [-] (-1,0) to (9,0);
\draw [-] (3,2) to (9,2);
\draw [-] (4,0) to node[scale=0.7] {$\overline{kp}+p$} (5,2);
\draw [-] (0,0) to node[scale=0.7] {$\overline{kp}$} (5,2);
\draw [-] (4,0) to node[scale=0.7] {$\overline{kp}+2p-N$} (9,2);
\draw [-] (8,0) to node[scale=0.7] {$\overline{kp}+3p-N$} (9,2);
\node at (0,-0.4) {$j$};
\node at (4,-0.4) {$j+1$};
\node at (8,-0.4) {$j+2$};
\end{tikzpicture}
\end{equation*}
In the first case mutation at $\overline{kp}+p$ gives $J_{\overline{kp}}$. Mutating the second case at $\overline{kp}+2p-N$ gives the first case so we can obtain $J_{\overline{kp}}$ between vertices $j$ and $j+2$ again. Since $\overline{kp}\equiv jp \mod q$ and $J_n$ is period $q$ we see that the arc obtained is indeed $J_{jp}$.

For $\tilde{J}$ we look at the leftmost arc at $q+i$ on the top boundary, labelled $\overline{kp}$. The situation is this:
\begin{equation}
\begin{tikzpicture}[every node/.style={fill=white}]
\draw[fill=black] (-2,0) circle [radius=0.1cm];
\draw[fill=black] (2,0) circle [radius=0.1cm];
\draw[fill=black] (4,0) circle [radius=0.1cm];
\draw[fill=black] (8,0) circle [radius=0.1cm];
\draw[fill=black] (9,2) circle [radius=0.1cm];
\draw[fill=black] (5,2) circle [radius=0.1cm];
\draw[fill=black] (-3,2) circle [radius=0.1cm];
\draw [-] (-3,0) to (9,0);
\draw [-] (-4,2) to (9,2);
\draw [-] (4,0) to node[scale=0.7] {$\overline{kp}+q-p$} (5,2);
\draw [-] (2,0) to (5,2);
\draw [-] (-2,0) to (5,2);
\draw [-] (-2,0) to (-3,2);
\draw [-] (8,0) to node[scale=0.7] {$\overline{kp}+q$} (5,2);
\draw [-] (8,0) to node[scale=0.7] {$\overline{kp}$} (9,2);
\node at (5,2.4) {$q+i-1$};
\node at (2,0.7) {$\ldots$};
\node at (0,0) {$\ldots$};
\node at (9,2.4) {$q+i$};
\end{tikzpicture}
\end{equation}
We mutate at each of the arcs at $q+i-1$ from left to right, except for the arc labelled $\overline{kp}+q$ which gives, near $q+i-1$,
\begin{equation}
\begin{tikzpicture}[every node/.style={fill=white}]
\draw[fill=black] (8,0) circle [radius=0.1cm];
\draw[fill=black] (9,2) circle [radius=0.1cm];
\draw[fill=black] (5,2) circle [radius=0.1cm];
\draw[fill=black] (-3,2) circle [radius=0.1cm];
\draw [-] (-3,0) to (9,0);
\draw [-] (-4,2) to (9,2);
\draw [-] (8,0) to node[scale=0.7] {$\overline{kp}+2q$} (-3,2);
\draw [-] (8,0) to node[scale=0.7] {$\overline{kp}+q$} (5,2);
\draw [-] (8,0) to node[scale=0.7] {$\overline{kp}$} (9,2);
\node at (5,2.4) {$q+i-1$};
\node at (-3,2.4) {$q+i-2$};
\node at (9,2.4) {$q+i$};
\node at (1,0.7) {$\ldots$};
\node at (3,0) {$\ldots$};
\end{tikzpicture}
\end{equation}
Mutation at $\overline{kp}+q$ will then give $\tilde{J}_{\overline{kp}}$. By a similar argument to the one at the start of this proof we have that $\tilde{J}_{\overline{kp}}=\tilde{J}_{-iq}$
\end{proof}
\end{lemma}
The $J$ and $\tilde{J}$ connect vertices on the same boundary by ``jumping" over one vertex. The arcs we haven't discussed so far are those jumping over more than one vertex. We define these as $J^{l-1}_{j}$, which starts at $j$ and jumps right over $l-1$ vertices to reach $j+l$, and $\tilde{J}_{i}^{m-1}$ which starts at $q+i$ and jumps left over $m-1$ vertices to reach $q+i-m$, as shown in Figure \ref{LongerJ}. Again we define these for $j=0,1,\ldots,q-1$ and $i=0,1,\ldots,p-1$.  We remark that $J^{q-1}_{j}$ and $\tilde{J}^{p-1}_{i}$ are the widest arcs, as $J^{q}_{j}$ and $\tilde{J}^{p}_{i}$ self-intersect.
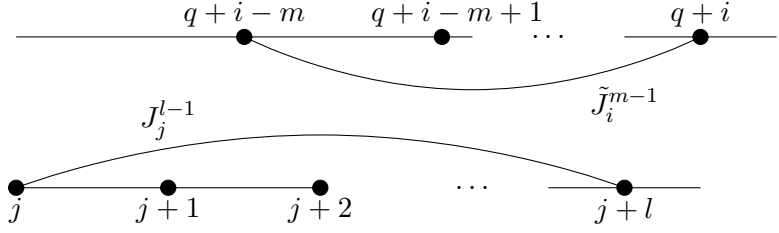
\begin{figure}
\centering
\begin{tikzpicture}
\draw[fill=black] (0,0) circle [radius=0.1cm];
\draw[fill=black] (2,0) circle [radius=0.1cm];
\draw[fill=black] (4,0) circle [radius=0.1cm];
\draw[fill=black] (3,2) circle [radius=0.1cm];
\draw[fill=black] (5.6,2) circle [radius=0.1cm];
\draw[fill=black] (9,2) circle [radius=0.1cm];
\draw[fill=black] (8,0) circle [radius=0.1cm];
\draw plot [smooth, tension=1] coordinates {(0,0) (4,0.7) (8,0)};
\draw plot [smooth, tension=1] coordinates {(3,2) (6,1.3) (9,2)};
\draw plot [smooth, tension=1] coordinates {(0,0) (2,0)};
\draw plot [smooth, tension=1] coordinates {(0,2) (6,2)};
\draw plot [smooth, tension=1] coordinates {(8,2) (10,2)};
\draw plot [smooth, tension=1] coordinates {(2,0) (4,0)};
\draw plot [smooth, tension=1] coordinates {(7,0) (9,0)};
\node at (6,0) {$\ldots$};
\node at (7,2) {$\ldots$};
\node at (2,0.9) {$J^{l-1}_{j}$};
\node at (8,1.1) {$\tilde{J}^{m-1}_{i}$};
\node at (0,-0.3) {$j$};
\node at (2,-0.3) {$j+1$};
\node at (4,-0.3) {$j+2$};
\node at (8,-0.3) {$j+l$};
\node at (3,2.3) {$q+i-m$};
\node at (5.8,2.3) {$q+i-m+1$};
\node at (9,2.3) {$q+i$};
\end{tikzpicture}\caption{The arcs $J^{l-1}_{j}$ and $\tilde{J}^{m-1}_{i}$.}\label{LongerJ}
\end{figure} 
For small $l$ we have 
\[
J^0_{j}=1, \qquad J^1_{j}=J_{jp}, \qquad \tilde{J}^0_{i}=1, \qquad \tilde{J}^1_{i}=\tilde{J}_{-iq}
\] 
and the following theorem allows us to calculate $J^l_{j}$ and $\tilde{J}^l_i$ for $l>1$.
\begin{theorem}
The arcs $J^l_j$ satisfy the recurrence relation
\[
J^{l-1}_{j}=J^{l-2}_{j}J_{(j+l-2)p}-J^{l-3}_{j}
\]
for $l=3,4,\ldots,q$, with initial values $J^l_0=1$ and $J^1_{j}=J_{jp}$. Similarly, the arcs $\tilde{J}^m_i$ satisfy
\[
\tilde{J}_i^{m-1}=\tilde{J}^{m-2}_i\tilde{J}_{-(i-m+2)q}-\tilde{J}^{m-3}_i
\] 
for $m=3,4,\ldots,p$, with initial values $\tilde{J}^0_{i}=1$ and $\tilde{J}^1_{i}=\tilde{J}_{-iq}$. Hence $J^l_{j}={D}^{l}_p(J_{jp})$, as defined in (\ref{determinant}), and they form a frieze (\ref{frieze}). Similarly we have $\tilde{J}^m_i=D^m_q(\tilde{J}_{-iq})$ forming a frieze (\ref{tildefrieze}).
\begin{proof}
We look at the quadrilateral with diagonals $J_{(j+l-2)p}$ and $J^{l-2}_{j}$:
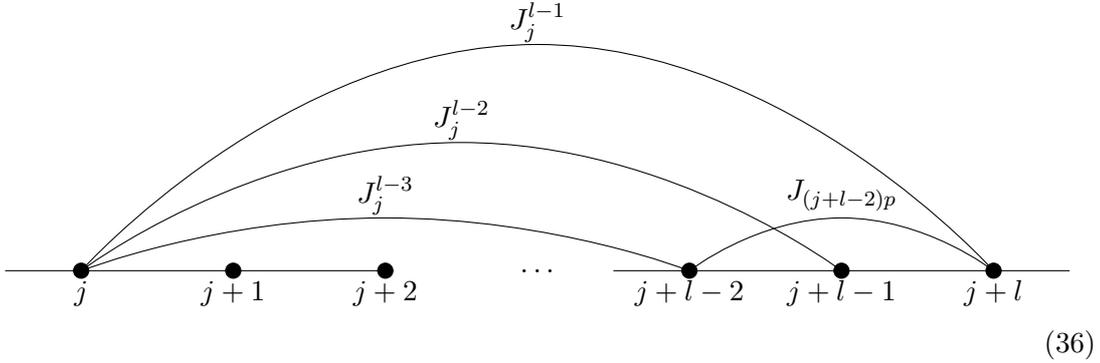
\begin{equation}\label{Ptolemypic}
\begin{tikzpicture}
\draw[fill=black] (0,0) circle [radius=0.1cm];
\draw[fill=black] (2,0) circle [radius=0.1cm];
\draw[fill=black] (4,0) circle [radius=0.1cm];
\draw[fill=black] (8,0) circle [radius=0.1cm];
\draw[fill=black] (10,0) circle [radius=0.1cm];
\draw[fill=black] (12,0) circle [radius=0.1cm];
\draw plot [smooth, tension=1] coordinates {(0,0) (6,3) (12,0)};
\draw plot [smooth, tension=1] coordinates {(0,0) (4,0.7) (8,0)};
\draw plot [smooth, tension=1] coordinates {(0,0) (5,1.7) (10,0)};
\draw plot [smooth, tension=1] coordinates {(8,0) (10,0.7) (12,0)};
\draw plot [smooth, tension=1] coordinates {(-1,0) (4,0)};
\draw plot [smooth, tension=1] coordinates {(7,0) (13,0)};
\node at (6,0) {$\ldots$};
\node at (6,3.3) {$J^{l-1}_{j}$};
\node at (10,1.0) {$J_{(j+l-2)p}$};
\node at (4,1.0) {$J^{l-3}_{j}$};
\node at (5,2.0) {$J^{l-2}_{j}$};
\node at (0,-0.3) {$j$};
\node at (2,-0.3) {$j+1$};
\node at (4,-0.3) {$j+2$};
\node at (8,-0.3) {$j+l-2$};
\node at (10,-0.3) {$j+l-1$};
\node at (12,-0.3) {$j+l$};
\end{tikzpicture}
\end{equation}
In this case the Ptolemy relation (\ref{ptolemy}) gives 
\[
J^{l-2}_{j}J_{(j+l-2)p}=J^{l-3}_{j}+J^{l-1}_{j}
\]
since the boundary arcs have the value $1$. The widest arcs are $J^{q-1}_{j}$, since the next arc, $J^q_{j}$, would self-intersect. The result for the $\tilde{J}^m_i$ arcs follows from a similar construction.  
\end{proof}
\end{theorem}
This proves the $\tilde{A}$ parts of Theorem \ref{findingallclustervarstheorem} and Proposition \ref{friezepropintro}. What remains is to prove the $\tilde{A}$ part of Proposition \ref{intropropclusterfrieze}.
\subsection{$\tilde{A}$ type cluster friezes}\label{Atypeclusterfriezesection}
Here we intend to show that the frieze pattern (\ref{friezeformulaintro}) and the friezes of Proposition \ref{friezepropintro} give a cluster frieze on the Auslander-Reiten quiver $\Gamma(\mathcal{C}_Q)$. They are already friezes, so we just need to ensure that they are connected by the relation (\ref{clusterfrieze}).

In this case there are two exceptional tubes given by $\lambda=0,1$. Firstly we need a description of the modules $B_{\lambda}$ and $B'_{\lambda}$, as given in \cite{assemdupont}. The vertex $e$ is required to be a sink so we first mutate the $\tilde{A}_{q,p}$ quiver at $0$ and take $e=0$. The portion of the quiver near $0$ now looks like
\begin{equation*}
\begin{tikzcd}
2q-N\arrow[r] & {q}\arrow[r] & {0} & p\arrow[l]
\end{tikzcd}
\end{equation*}
with cluster variables
\begin{equation}\label{newnearzero}
\begin{tikzcd}
x_{2q-N}\arrow[r] & x_{q}\arrow[r] & x_{N} & x_p\arrow[l]
\end{tikzcd}
\end{equation}
and \cite{assemdupont} gives $B_0=P_{q}$ and  $B'_0=P_{p}[1]$. By definition $X_{P_{p}[1]}=x_p$. $X_{P_0}=X_{S_0}$ has $x_N$ as its denominator, hence is given by mutating (\ref{newnearzero}) at $N$, so $X_{P_0}=x_0$. Similarly $X_{P_q}$ has denominator $x_qx_N$ so is given by performing $\mu_q\mu_{N}$ on (\ref{newnearzero}), so $X_{P_q}=x_{-p}$. Collecting this we have that (\ref{clusterfrieze}) is
\[
X_{N_0}=\frac{X_{P_{q}}+X_{P_{p}[1]}}{X_{P_{N}}}=\frac{x_{-p}+x_{p}}{x_{0}}=J_{0}.
\]
The map $M\mapsto M[1]$ gives an automorphism of the cluster algebra, so
\[
X_{N_0[j]}=\frac{X_{P_{q}[j]}+X_{P_{p}[j+1]}}{X_{P_{N}[j]}}=\frac{x_{-p+jN}+x_{p+jN}}{x_{-1+jN}}=J_{jN}=J_{jp}.
\]
Conversely the modules $B_1=P_{p}$ and $B'_1=P_{q}[1]$ give
\[
X_{N_1[j]}=\frac{X_{P_{q}[j]}+X_{P_{p}[j+1]}}{X_{P_{N}[j]}}=\frac{x_{-q+jN}+x_{q+jN}}{x_{jN}}=\tilde{J}_{jN}=\tilde{J}_{jq}.
\]
This proves Proposition \ref{intropropclusterfrieze} in that $\tilde{A}$ case.
\section{Triangulated surfaces and the cluster frieze in $\tilde{D}$ type}\label{Dtypesection}
In this section we look at the construction of $\tilde{D}$ quivers as triangulations of discs with two punctures. We first show how to construct the arcs corresponding to the cluster map variables $X^i_n$. We then construct a periodic frieze of cluster variables, with the first row given by the periodic quantities $J'_n$, as in the $\tilde{A}$ case. We show that there are only $3$ exceptional arcs outside of these. Finally we identify the two friezes constructed here as the cluster friezes given in \cite{assemdupont}. 
\subsection{$\tilde{D}$ type cluster algebras as triangulated surfaces}
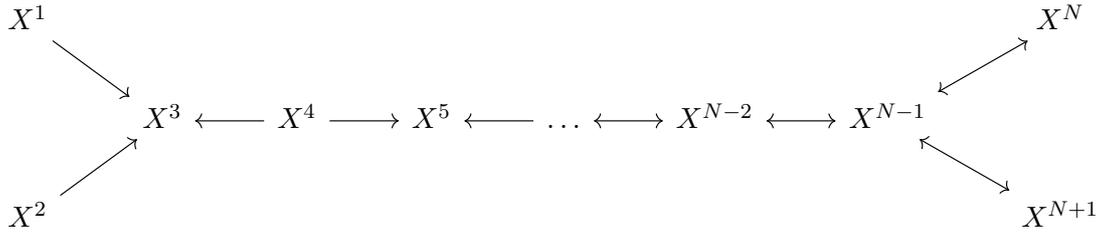
\begin{figure}
\centering
\begin{equation*}
\begin{tikzcd}
X^1\arrow[dr] & & & & & & & X^N\\
& X^3 & X^4\arrow[l]\arrow[r] & X^5 & \arrow[l]\ldots\arrow[r, leftrightarrow] & X^{N-2}\arrow[r, leftrightarrow] & X^{N-1}\arrow[ur, leftrightarrow]\arrow[dr, leftrightarrow] \\
X^2\arrow[ur] & & & & & & & X^{N+1}
\end{tikzcd}
\end{equation*}
\caption{The $\tilde{D}_N$ quiver.}\label{Dquiver}
\end{figure}

As discussed in Subsection \ref{dynamicalsystemssubsection} we take the bipartite orientation of the $\tilde{D}$ diagram as shown in Figure \ref{Dquiver}. To ensure that this is bipartite the orientation of the arrows at the right end of the diagram depend on the parity of $N$, which we have signified with double ended arrows.
The cluster map 
\begin{equation}\label{Dtypeclustermap}
\varphi:
\begin{pmatrix}
X^1_n, & X^2_n, &\ldots &X^{N+1}_n
\end{pmatrix}
\mapsto
\begin{pmatrix}
X^1_{n+1}, & X^2_{n+1}, &\ldots &X^{N+1}_{n+1}
\end{pmatrix}
\end{equation}
is obtained by applying
\[
\mu:=\mu_{\mathrm{source}}\circ\mu_{\mathrm{sink}}
\]
to $Q$, where
\[
\mu_{\mathrm{source}}:=\mu_1\mu_2\mu_4\ldots, \qquad \mu_{\mathrm{sink}}:=\mu_3\mu_5\mu_7\ldots
\]
which are the compositions of mutations at the sources and sinks in Figure \ref{Dquiver}. As shown in \cite{clustersandtriangulated1} we can see this quiver as a triangulation of a disk with $2$ marked points inside and $N-2$ marked points on the boundary. An example of this is seen in Figure \ref{Dtypebefore}. 

In order to determine all of the arcs that can occur in a triangulation, we firstly note that every arc $\gamma$ connecting boundary vertices splits the disc into two connected components. There are two possibilities: either the punctures $p_1$ and $p_2$ lie in the same connected component or in different components. We denote these sets of arcs by:
\[
\Gamma_{2,0}:=\{\gamma \mid \gamma \textrm{ connects boundary vertices and } p_1 \textrm{ and } p_2 \textrm{ appear on the same side of } \gamma\}
\]
\[
\Gamma_{1,1}:=\{\gamma \mid \gamma \textrm{ connects boundary vertices and } p_1 \textrm{ and } p_2 \textrm{ appear on different sides of } \gamma\}
\]
There are also arcs connecting the boundary to the punctures, which we call
\[
\Gamma_{\mathrm{punc}}:=\{\gamma \mid \gamma \textrm{ connects a puncture to the boundary}\}
\]
Outside of these three sets there are three exceptional arcs involving only the punctures (\ref{3exceptionalarcs}), which we call $\Gamma_{\mathrm{except}}$.
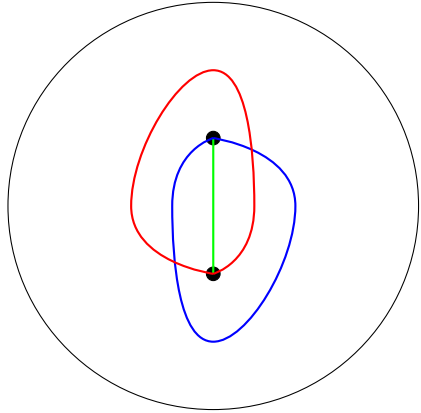
\begin{equation}\label{3exceptionalarcs}
\begin{tikzpicture}[scale=0.9, every node/.style={fill=white}]
\draw (0,0) circle [radius=3.0];
\draw[fill=black] (0,1) circle [radius=0.1cm];
\draw[fill=black] (0,-1) circle [radius=0.1cm];
\draw[green, thick] plot [smooth, tension=1] coordinates { (0,1) (0,-1) };
\draw[blue, thick] plot [smooth, tension=1] coordinates { (0,1) (-0.6,0) (0,-2) (1.2,0) (0,1) };
\draw[red, thick] plot [smooth, tension=1] coordinates {(0,-1) (-1.2,0) (0,2) (0.6,0) (0,-1)};
\end{tikzpicture}
\end{equation}
\begin{lemma}\label{Dsurfacepossiblearcs}
The arcs of $\Gamma_{1,1}$ look like those shown in Figure \ref{just2boundaryverticesfinal}. They start from a vertex $v_i$ and encircle the line $L$ and the punctures $m\in \mathbb{Z}$ times, as shown in Figure \ref{just2boundaryvertices} in blue for $m>0$. After this the curve crosses $L$ and takes the only possible path to $v_j$ (red). For $m<0$ the curve is shown by reflecting this picture in the vertical axis. We denote these by $\gamma(v_i,v_j,m)$.  

The arcs of $\Gamma_{\mathrm{punc}}$ are the arcs that form self-folded triangles with a $\gamma(v_i,v_i,m)$.
\begin{proof}
We firstly look at an arc $\gamma\in \Gamma_{1,1}$ connecting two boundary vertices $v_i$ and $v_j$. Since this arc separates the two punctures it necessarily crosses the line joining them, which we call $L$. The portion of $\gamma$ before this crossing can only live in the annulus obtained by removing an area enclosing the two punctures and $L$, as shown in Figure \ref{just2boundaryvertices} on the left, so it simply orbits the drawn ellipse $m\in \mathbb{Z}$ times, where we take $m>0$ to mean that the blue portion of the arc travels anticlockwise, starting from $v_i$, while $m<0$ means clockwise. 

When $\gamma$ finally crosses $L$, as shown on the right of Figure \ref{just2boundaryvertices} (we draw this new portion of $\gamma$ in red for clarity), there is only one choice, to follow the maze back to $v_j$. The only remaining issue is that perhaps the curve crosses $L$ from the other direction, i.e. the curve of Figure \ref{just2boundaryvertices} is blue to the right of $L$ and red to the left. In this case we just need to permute $i\leftrightarrow j$ to obtain a curve of the right form. 

Finally we note that every arc of $\Gamma_{\mathrm{punc}}$ appears as the internal arc of a self-folded triangle. The other arc of this self-folded triangle, as we have just shown, is of the form $\gamma(v_i,v_i,m)$.
\begin{figure}
\centering
\begin{tikzpicture}[scale=0.85, every node/.style={fill=white}]
\draw (0,0) circle [radius=4.0];
\draw[fill=black] (0,1) circle [radius=0.1cm];
\draw[fill=black] (0,-1) circle [radius=0.1cm];
\draw[fill=black] (0,-4) circle [radius=0.1cm];
\draw[fill=black] (-3.46,2) circle [radius=0.1cm];
\node[scale=1] at (-4,2){$v_j$};
\node[scale=1] at (0,-4.4){$v_i$};
\node[scale=0.7] at (0.2,0){$L$};
\draw[dotted] (0,-1) to (0,1);
\draw (0,0) ellipse (1cm and 2cm);
\draw[blue, thick] plot [smooth, tension=1] coordinates {(0,-4) (3.5,0) (0,3.5) (-3.5,0) (0,-3) (2.9,0) (0,2.9) (-2.9,0) (0,-2.5) (2.5,0) (0,2.5)};
\end{tikzpicture}
\qquad
\begin{tikzpicture}[every node/.style={fill=white}]
\draw[fill=black] (0,1) circle [radius=0.1cm];
\draw[fill=black] (0,-1) circle [radius=0.1cm];
\node[scale=0.7] at (0.2,0.65){$L$};
\draw[dotted] (0,-1) to (0,1);
\draw[blue, thick] plot [smooth, tension=1] coordinates { (-1,3.5) (-3.5,0) (0,-3) (2.9,0) (0,2.9) (-2.9,0) (0,-2.5) (2.5,0) (0,2.3) (-1.8,0.8) (0,0)};
\draw[red, thick] plot [smooth, tension=1] coordinates {(0,0) (1,-0.8) (0,-2) (-2,0)};
\end{tikzpicture}
\caption{The portion of an arc (blue) before crossing the line $L$ and the behaviour just after crossing (red).}
\label{just2boundaryvertices}
\end{figure}
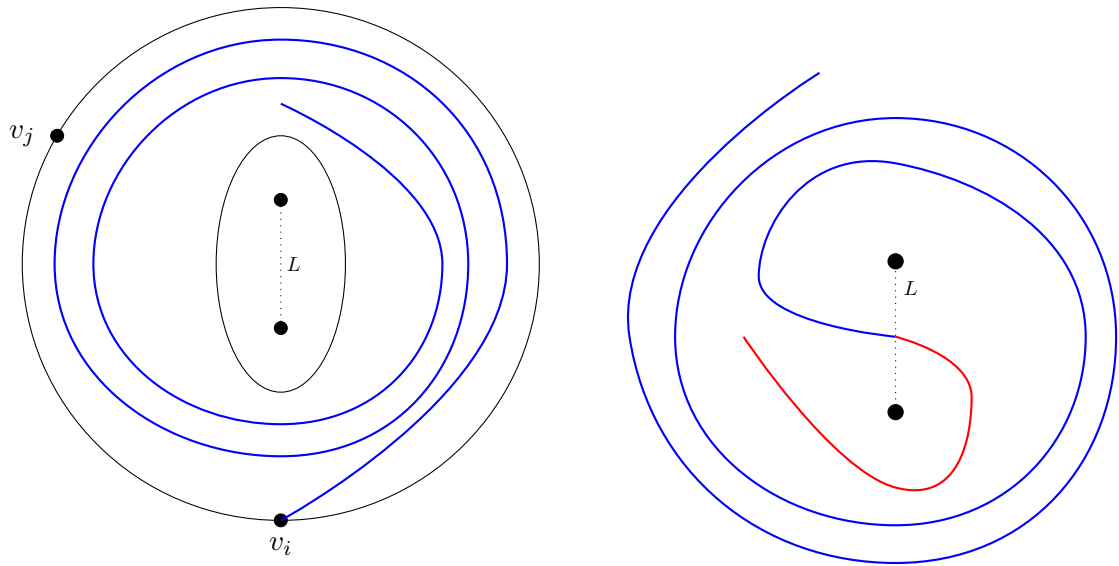
\begin{figure}
\centering
\begin{tikzpicture}[every node/.style={fill=white}]
\draw (0,0) circle [radius=4.0];
\draw[fill=black] (0,1) circle [radius=0.1cm];
\draw[fill=black] (0,-1) circle [radius=0.1cm];
\draw[fill=black] (0,-4) circle [radius=0.1cm];
\draw[fill=black] (-3.46,2) circle [radius=0.1cm];
\node[scale=1] at (-3.9,2){$v_j$};
\node[scale=1] at (0,-4.4){$v_i$};
\node[scale=0.7] at (0.2,0.65){$L$};
\draw[dotted] (0,-1) to (0,1);
\draw[blue, thick] plot [smooth, tension=1] coordinates {(0,-4) (3.5,0) (0,3.5) (-3.5,0) (0,-3) (2.9,0) (0,2.9) (-2.9,0) (0,-2.5) (2.5,0) (0,2.3) (-1.8,0.8) (0,0)};
\draw[red, thick] plot [smooth, tension=1] coordinates {(0,0) (1,-0.8) (0,-2) (-2,0) (-1,2.4) (1,2.2) (2.6,0.7) (2.0,-1.7) (0,-2.7) (-2,-2) (-3.2,0) (-2,2.6) (0,3.2) (2,2.6) (3.2,0) (2,-2.4) (0,-3.3) (-2.1,-2.9) (-3.7,0) (-3.3,1.6) (-3.46,2) };
\end{tikzpicture}
\caption{The full curve after the red portion takes the only path to reach $v_j$.}
\label{just2boundaryverticesfinal}
\end{figure}
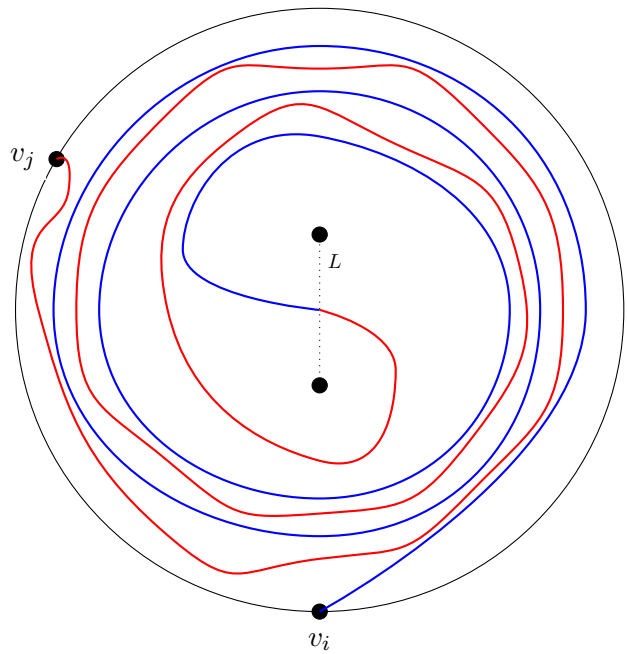
\end{proof}
\end{lemma}
Next we try to identify the arcs obtained by the cluster map. Our initial triangulation is shown in Figure \ref{Dtypebefore} for $N=8$. Its construction can be tentatively described in two steps: 
\begin{enumerate}[(i)]
\item Draw two self folded triangles; one at $v_1$ and $p_1$ and one at $v_6$ and $p_2$.
\item Draw arcs $i$ from $v_{i-2}$ to $v_{i-1}$ for $i=3,4,\ldots, N-1$.
\end{enumerate}
We remark that in the second step there is only one choice for each of the arcs. We need to see how this triangulation changes under $\mu$, which we perform in several steps, as shown in Figure \ref{fullDmutation}. 
\begin{figure}
  \begin{subfigure}{0.31\textwidth}
 \begin{tikzpicture}[scale=0.85, every node/.style={fill=white}]
\draw (0,0) circle [radius=4.0];
\draw[fill=black] (0,4) circle [radius=0.1cm];
\draw[fill=black] (0,2) circle [radius=0.1cm];
\draw[fill=black] (0,-2) circle [radius=0.1cm];
\draw[fill=black] (0,-4) circle [radius=0.1cm];
\draw[fill=black] (3.46,2) circle [radius=0.1cm];
\draw[fill=black] (-3.46,2) circle [radius=0.1cm];
\draw[fill=black] (3.46,-2) circle [radius=0.1cm];
\draw[fill=black] (-3.46,-2) circle [radius=0.1cm];
\draw plot [smooth, tension=1] coordinates {(0,-4) (0,-2)};
\draw plot [smooth, tension=1] coordinates {(0,4) (0,2)};
\draw plot [smooth, tension=1] coordinates {(0,4) (-0.7,3) (0,1.6) (0.7,3) (0,4)};
\draw plot [smooth, tension=1] coordinates {(0,-4) (-0.7,-3) (0,-1.6) (0.7,-3) (0,-4)};
\draw plot [smooth, tension=0.9] coordinates {(0,-4) (-1,-3) (0,-1.4) (3.46,-2)};
\draw plot [smooth, tension=0.9] coordinates {(0,4) (1,3) (0,1.4) (-3.46,2)};
\draw plot [smooth, tension=0.9] coordinates {(3.46,2) (0,1) (-3.46,2)};
\draw plot [smooth, tension=0.9] coordinates {(3.46,-2) (0,-1) (-3.46,-2)};
\draw plot [smooth, tension=0.9] coordinates {(3.46,2) (0,0) (-3.46,-2)};
\node[scale=0.7] at (0,-3){1};
\node[scale=0.7] at (0.7,-3){2};
\node[scale=0.7] at (-1,-2.7){3};
\node[scale=0.7] at (0,-1){4};
\node[scale=0.7] at (0,0){5};
\node[scale=0.7] at (0,1){6};
\node[scale=0.7] at (1,2.7){7};
\node[scale=0.7] at (-0.7,3){8};
\node[scale=0.7] at (0,3){9};
\node[scale=0.7] at (-3.9,-2){$v_3$};
\node[scale=0.7] at (-3.9,2){$v_5$};
\node[scale=0.7] at (3.9,2){$v_4$};
\node[scale=0.7] at (0,-4.4){$v_1$};
\node[scale=0.7] at (0,4.4){$v_6$};
\node[scale=0.7] at (3.9,-2){$v_2$};
\node[scale=0.7] at (-0.2,-2.35){$p_1$};
\node[scale=0.7] at (0.2,2.35){$p_2$};
\end{tikzpicture}
    \caption{The initial triangulation for $N=8$.} \label{Dtypebefore}
  \end{subfigure}%
  \hspace{40mm}  
  \begin{subfigure}{0.31\textwidth}
   \begin{tikzpicture}[scale=0.85, every node/.style={fill=white}]
\draw (0,0) circle [radius=4.0];
\draw[fill=black] (0,4) circle [radius=0.1cm];
\draw[fill=black] (0,2) circle [radius=0.1cm];
\draw[fill=black] (0,-2) circle [radius=0.1cm];
\draw[fill=black] (0,-4) circle [radius=0.1cm];
\draw[fill=black] (3.46,2) circle [radius=0.1cm];
\draw[fill=black] (-3.46,2) circle [radius=0.1cm];
\draw[fill=black] (3.46,-2) circle [radius=0.1cm];
\draw[fill=black] (-3.46,-2) circle [radius=0.1cm];
\draw plot [smooth, tension=1] coordinates {(0,-4) (0,-2)};
\draw plot [smooth, tension=1] coordinates {(0,4) (0,2)};
\draw plot [smooth, tension=1] coordinates {(0,4) (-0.7,3) (0,1.6) (0.7,3) (0,4)};
\draw plot [smooth, tension=1] coordinates {(0,-4) (-0.7,-3) (0,-1.6) (0.7,-3) (0,-4)};
\draw plot [smooth, tension=0.9] coordinates {(0,-4) (1,-3) (0,-1.4) (-3.46,-2)};
\draw plot [smooth, tension=0.9] coordinates {(0,4) (-1,3) (0,1.4) (3.46,2)};
\draw plot [smooth, tension=0.9] coordinates {(3.46,2) (0,1) (-3.46,2)};
\draw plot [smooth, tension=0.9] coordinates {(3.46,-2) (0,-1) (-3.46,-2)};
\draw plot [smooth, tension=0.9] coordinates {(-3.46,2) (0,0) (3.46,-2)};
\node[scale=0.7] at (0,-3){1};
\node[scale=0.7] at (-0.7,-3){2};
\node[scale=0.7] at (1,-2.7){3};
\node[scale=0.7] at (0,-1){4};
\node[scale=0.7] at (0,0){5};
\node[scale=0.7] at (0,1){6};
\node[scale=0.7] at (-1,2.7){7};
\node[scale=0.7] at (0.7,3){8};
\node[scale=0.7] at (0,3){9};
\end{tikzpicture}
    \caption{After performing $\mu_7\mu_5\mu_3$.} \label{fig:1b}
  \end{subfigure}%

  \begin{subfigure}{0.31\textwidth}
    \begin{tikzpicture}[scale=0.85, every node/.style={fill=white}]
\draw (0,0) circle [radius=4.0];
\draw[fill=black] (0,4) circle [radius=0.1cm];
\draw[fill=black] (0,2) circle [radius=0.1cm];
\draw[fill=black] (0,-2) circle [radius=0.1cm];
\draw[fill=black] (0,-4) circle [radius=0.1cm];
\draw[fill=black] (3.46,2) circle [radius=0.1cm];
\draw[fill=black] (-3.46,2) circle [radius=0.1cm];
\draw[fill=black] (3.46,-2) circle [radius=0.1cm];
\draw[fill=black] (-3.46,-2) circle [radius=0.1cm];
\draw plot [smooth, tension=1] coordinates {(0,-4) (0,-2)};
\draw plot [smooth, tension=1] coordinates {(0,4) (0,2)};
\draw plot [smooth, tension=1] coordinates {(0,4) (-0.7,3) (0,1.6) (0.7,3) (0,4)};
\draw plot [smooth, tension=1] coordinates {(0,-2) (-3.46,-2)};
\draw plot [smooth, tension=0.9] coordinates {(0,-4) (1,-3) (0,-1.4) (-3.46,-2)};
\draw plot [smooth, tension=0.9] coordinates {(0,4) (-1,3) (0,1.4) (3.46,2)};
\draw plot [smooth, tension=0.9] coordinates {(3.46,2) (0,1) (-3.46,2)};
\draw plot [smooth, tension=0.9] coordinates {(3.46,-2) (0,-1) (-3.46,-2)};
\draw plot [smooth, tension=0.9] coordinates {(-3.46,2) (0,0) (3.46,-2)};
\node[scale=0.7] at (0,-3){1};
\node[scale=0.7] at (-0.7,-2){2};
\node[scale=0.7] at (1,-2.7){3};
\node[scale=0.7] at (0,-1){4};
\node[scale=0.7] at (0,0){5};
\node[scale=0.7] at (0,1){6};
\node[scale=0.7] at (-1,2.7){7};
\node[scale=0.7] at (0.7,3){8};
\node[scale=0.7] at (0,3){9};
\end{tikzpicture}
    \caption{After $\mu_2$.} \label{Dtypenext}
  \end{subfigure}
\hspace{40mm} 
\begin{subfigure}{0.31\textwidth}
 \begin{tikzpicture}[scale=0.85, every node/.style={fill=white}]
\draw (0,0) circle [radius=4.0];
\draw[fill=black] (0,4) circle [radius=0.1cm];
\draw[fill=black] (0,2) circle [radius=0.1cm];
\draw[fill=black] (0,-2) circle [radius=0.1cm];
\draw[fill=black] (0,-4) circle [radius=0.1cm];
\draw[fill=black] (3.46,2) circle [radius=0.1cm];
\draw[fill=black] (-3.46,2) circle [radius=0.1cm];
\draw[fill=black] (3.46,-2) circle [radius=0.1cm];
\draw[fill=black] (-3.46,-2) circle [radius=0.1cm];
\draw plot [smooth, tension=1] coordinates {(-3.46,-2) (-0.5,-1.7) (0.5,-2) (-0.5,-2.3) (-3.46,-2)};
\draw plot [smooth, tension=1] coordinates {(0,4) (0,2)};
\draw plot [smooth, tension=1] coordinates {(0,4) (-0.7,3) (0,1.6) (0.7,3) (0,4)};
\draw plot [smooth, tension=1] coordinates {(0,-2) (-3.46,-2)};
\draw plot [smooth, tension=0.9] coordinates {(0,-4) (1,-3) (0,-1.4) (-3.46,-2)};
\draw plot [smooth, tension=0.9] coordinates {(0,4) (-1,3) (0,1.4) (3.46,2)};
\draw plot [smooth, tension=0.9] coordinates {(3.46,2) (0,1) (-3.46,2)};
\draw plot [smooth, tension=0.9] coordinates {(3.46,-2) (0,-1) (-3.46,-2)};
\draw plot [smooth, tension=0.9] coordinates {(-3.46,2) (0,0) (3.46,-2)};
\node[scale=0.7] at (0.5,-2){1};
\node[scale=0.7] at (-0.7,-2){2};
\node[scale=0.7] at (1,-2.7){3};
\node[scale=0.7] at (0,-1){4};
\node[scale=0.7] at (0,0){5};
\node[scale=0.7] at (0,1){6};
\node[scale=0.7] at (-1,2.7){7};
\node[scale=0.7] at (0.7,3){8};
\node[scale=0.7] at (0,3){9};
\end{tikzpicture}
    \caption{After $\mu_1$.}
  \end{subfigure}

\begin{subfigure}{0.31\textwidth}
\begin{tikzpicture}[scale=0.85, every node/.style={fill=white}]
\draw (0,0) circle [radius=4.0];
\draw[fill=black] (0,4) circle [radius=0.1cm];
\draw[fill=black] (0,2) circle [radius=0.1cm];
\draw[fill=black] (0,-2) circle [radius=0.1cm];
\draw[fill=black] (0,-4) circle [radius=0.1cm];
\draw[fill=black] (3.46,2) circle [radius=0.1cm];
\draw[fill=black] (-3.46,2) circle [radius=0.1cm];
\draw[fill=black] (3.46,-2) circle [radius=0.1cm];
\draw[fill=black] (-3.46,-2) circle [radius=0.1cm];
\draw plot [smooth, tension=1] coordinates {(-3.46,-2) (-0.5,-1.7) (0.5,-2) (-0.5,-2.3) (-3.46,-2)};
\draw plot [smooth, tension=1] coordinates {(0,4) (0,2)};
\draw plot [smooth, tension=1] coordinates {(0,4) (-0.7,3) (0,1.6) (0.7,3) (0,4)};
\draw plot [smooth, tension=1] coordinates {(0,-2) (-3.46,-2)};
\draw plot [smooth, tension=0.9] coordinates {(0,-4) (1,-3) (0,-1.4) (-3.46,-2)};
\draw plot [smooth, tension=0.9] coordinates {(0,4) (-1,3) (0,1.4) (3.46,2)};
\draw plot [smooth, tension=0.9] coordinates {(0,4) (-2,3) (0,0.7) (3.46,-2)};
\draw plot [smooth, tension=0.9] coordinates {(0,-4) (2,-3) (0,-0.7) (-3.46,2)};
\draw plot [smooth, tension=0.9] coordinates {(-3.46,2) (0,0) (3.46,-2)};
\node[scale=0.7] at (0.5,-2){1};
\node[scale=0.7] at (-0.7,-2){2};
\node[scale=0.7] at (1,-2.7){3};
\node[scale=0.7] at (0,-0.7){4};
\node[scale=0.7] at (0,0){5};
\node[scale=0.7] at (0,0.7){6};
\node[scale=0.7] at (-1,2.7){7};
\node[scale=0.7] at (0.7,3){8};
\node[scale=0.7] at (0,3){9};
\end{tikzpicture}
\caption{After $\mu_6\mu_4$.} \label{Dtypenext2}
  \end{subfigure}
\hspace{40mm} 
\begin{subfigure}{0.31\textwidth}
\begin{tikzpicture}[scale=0.85, every node/.style={fill=white}]
\draw (0,0) circle [radius=4.0];
\draw[fill=black] (0,4) circle [radius=0.1cm];
\draw[fill=black] (0,2) circle [radius=0.1cm];
\draw[fill=black] (0,-2) circle [radius=0.1cm];
\draw[fill=black] (0,-4) circle [radius=0.1cm];
\draw[fill=black] (3.46,2) circle [radius=0.1cm];
\draw[fill=black] (-3.46,2) circle [radius=0.1cm];
\draw[fill=black] (3.46,-2) circle [radius=0.1cm];
\draw[fill=black] (-3.46,-2) circle [radius=0.1cm];
\draw plot [smooth, tension=1] coordinates {(-3.46,-2) (-0.5,-1.7) (0.5,-2) (-0.5,-2.3) (-3.46,-2)};
\draw plot [smooth, tension=1] coordinates {(3.46,2) (0.5,1.7) (-0.5,2) (0.5,2.3) (3.46,2)};
\draw plot [smooth, tension=1] coordinates {(0,-2) (-3.46,-2)};
\draw plot [smooth, tension=1] coordinates {(0,2) (3.46,2)};
\draw plot [smooth, tension=0.9] coordinates {(0,-4) (1,-3) (0,-1.4) (-3.46,-2)};
\draw plot [smooth, tension=0.9] coordinates {(0,4) (-1,3) (0,1.4) (3.46,2)};
\draw plot [smooth, tension=0.9] coordinates {(0,4) (-2,3) (0,0.7) (3.46,-2)};
\draw plot [smooth, tension=0.9] coordinates {(0,-4) (2,-3) (0,-0.7) (-3.46,2)};
\draw plot [smooth, tension=0.9] coordinates {(-3.46,2) (0,0) (3.46,-2)};
\node[scale=0.7] at (0.5,-2){1};
\node[scale=0.7] at (-0.7,-2){2};
\node[scale=0.7] at (1,-2.7){3};
\node[scale=0.7] at (0,-0.7){4};
\node[scale=0.7] at (0,0){5};
\node[scale=0.7] at (0,0.7){6};
\node[scale=0.7] at (-1,2.7){7};
\node[scale=0.7] at (0.7,2){8};
\node[scale=0.7] at (-0.5,2){9};
\node[scale=0.7] at (-3.9,-2){$v_1$};
\node[scale=0.7] at (-3.9,2){$v_3$};
\node[scale=0.7] at (3.9,2){$v_6$};
\node[scale=0.7] at (0,-4.4){$v_2$};
\node[scale=0.7] at (0,-2.35){$p_1$};
\node[scale=0.7] at (0,2.35){$p_2$};
\node[scale=0.7] at (0,4.4){$v_5$};
\node[scale=0.7] at (3.9,-2){$v_4$};
\end{tikzpicture}
\caption{After $\mu_9\mu_8$.} \label{laststep}
  \end{subfigure}
\caption{A full application of $\mu$ to the initial quiver.} \label{fullDmutation}
\end{figure}
In Figure \ref{laststep}, after a complete application of $\mu$, we can see that the boundary vertices of the self-folded triangles have moved clockwise. We have also rotated the boundary vertices to obtain a new labelling for the boundary of $\mu(Q)$, so now the new quiver can be described by the same two steps as before.

In Figure \ref{incompleteDtype} we display the full triangulation for $\mu^2(Q)$ and some of the triangulation for $\mu^3(Q)$ (because the full picture would be too busy to be helpful) which are also constructed by the same two steps, if we relabel as before.
\begin{figure}
\begin{tikzpicture}[every node/.style={fill=white}]
\draw (0,0) circle [radius=4.0];
\draw[fill=black] (0,4) circle [radius=0.1cm];
\draw[fill=black] (0,2) circle [radius=0.1cm];
\draw[fill=black] (0,-2) circle [radius=0.1cm];
\draw[fill=black] (0,-4) circle [radius=0.1cm];
\draw[fill=black] (3.46,2) circle [radius=0.1cm];
\draw[fill=black] (-3.46,2) circle [radius=0.1cm];
\draw[fill=black] (3.46,-2) circle [radius=0.1cm];
\draw[fill=black] (-3.46,-2) circle [radius=0.1cm];
\draw plot [smooth, tension=1] coordinates {(-3.46,2) (0,-2)};
\draw[red, thick] plot [smooth, tension=1] coordinates {(-3.46,2) (-0.1,-1.6) (-0.3,-1.9) (-3.46,2)};
\draw[green, thick] plot [smooth, tension=1] coordinates {(-3.46,2) (0.3,-1.7) (-0.3,-2.6) (-2,-1.0) (-3.46,-2)};
\draw[green, thick] plot [smooth, tension=1] coordinates {(3.46,-2) (-0.3,1.7) (0.3,2.6) (2,1.0) (3.46,2)};
\draw[orange, thick] plot [smooth, tension=1] coordinates {(3.46,2) (-0.8,2.45) (2.4,-2.1) (0,-4)};
\draw[orange, thick] plot [smooth, tension=1] coordinates {(-3.46,-2) (0.8,-2.45) (-2.4,2.1) (0,4)};
\draw[blue, thick] plot [smooth, tension=1] coordinates {(0,-4) (2,-2) (0,0)};
\draw[blue, thick] plot[smooth, tension=1] coordinates {(0,4) (-2,2) (0,0)};
\draw plot [smooth, tension=1] coordinates {(3.46,-2) (0,2)};
\draw[red, thick] plot [smooth, tension=1] coordinates {(3.46,-2) (0.1,1.6) (0.3,1.9) (3.46,-2)};
\end{tikzpicture}
\qquad
\begin{tikzpicture}[every node/.style={fill=white}]
\draw (0,0) circle [radius=4.0];
\draw[fill=black] (0,4) circle [radius=0.1cm];
\draw[fill=black] (0,2) circle [radius=0.1cm];
\draw[fill=black] (0,-2) circle [radius=0.1cm];
\draw[fill=black] (0,-4) circle [radius=0.1cm];
\draw[fill=black] (3.46,2) circle [radius=0.1cm];
\draw[fill=black] (-3.46,2) circle [radius=0.1cm];
\draw[fill=black] (3.46,-2) circle [radius=0.1cm];
\draw[fill=black] (-3.46,-2) circle [radius=0.1cm];
\draw plot [smooth, tension=1] coordinates {(0,4) (-2.5,1) (0,-2)};
\draw plot [smooth, tension=1] coordinates {(0,-4) (2.5,-1) (0,2)};
\draw[red, thick] plot [smooth, tension=1] coordinates {(0,4) (-2.3,1) (0.4,-2.1) (-0.1,-2.25) };
\draw[red, thick] plot[smooth, tension=1] coordinates {(-0.1,-2.25) (-2.6,1) (0,4)};
\draw[green, thick] plot [smooth, tension=1] coordinates {(0,4) (-2.1,1) (0.6,-1.9) (0.6,-2.7) };
\draw[green, thick] plot [smooth, tension=1] coordinates {(0.6,-2.7) (-1,-1.8) (-2.8,0) (-3.46,2)};
\draw[red, thick] plot [smooth, tension=1] coordinates {(0,-4) (2.3,-1) (-0.4,2.1) (0.1,2.25) };
\draw[red, thick] plot[smooth, tension=1] coordinates {(0.1,2.25) (2.6,-1) (0,-4)};
\draw[orange, thick] plot [smooth, tension=1] coordinates {(-3.46,2) (-2.5,-2) (0,-3.15) (1.3,-2.3) (0,0)};
\draw[orange, thick] plot [smooth, tension=1] coordinates {(3.46,2) (0.5,3.15) (-1.3,2.3) (0,0)};
\draw[blue, thick] plot [smooth, tension=1] coordinates {(3.46,2) (0.7,2.9) (-0.9,2.3) };
\end{tikzpicture}
\caption{Further applications: $\mu^2(Q)$ on the left and an incomplete drawing of $\mu^3(Q)$ on the right.} \label{incompleteDtype}
\end{figure}
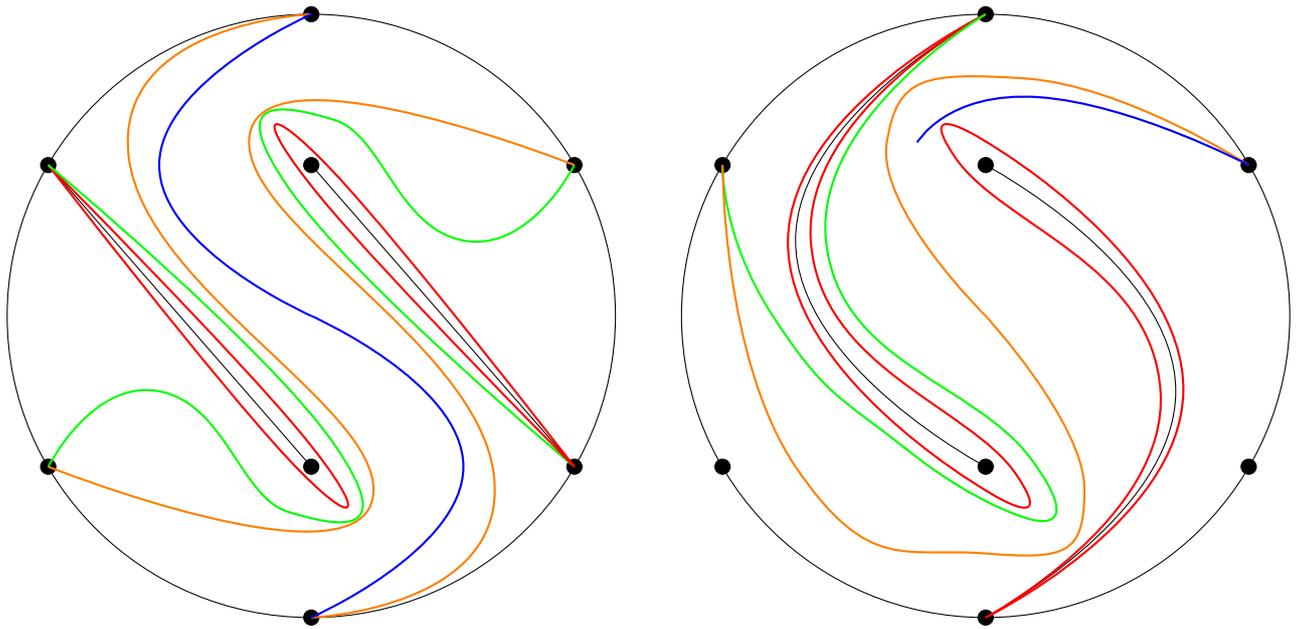
\begin{figure}
\centering
\begin{tikzpicture}[every node/.style={fill=white}]
\draw (0,0) circle [radius=4.0];
\draw[fill=black] (0,4) circle [radius=0.1cm];
\draw[fill=black] (0,2) circle [radius=0.1cm];
\draw[fill=black] (0,-2) circle [radius=0.1cm];
\draw[fill=black] (0,-4) circle [radius=0.1cm];
\draw[fill=black] (3.46,2) circle [radius=0.1cm];
\draw[fill=black] (-3.46,2) circle [radius=0.1cm];
\draw[fill=black] (3.46,-2) circle [radius=0.1cm];
\draw[fill=black] (-3.46,-2) circle [radius=0.1cm];
\draw plot [smooth, tension=1] coordinates {(0,-4) (3,1.8) (0,3)};
\draw plot [smooth, tension=1] coordinates {(0,3) (-1.3,1.2) (0,-2)};
\draw plot [smooth, tension=1] coordinates {(0,4) (-3,-1.8) (0,-3)};
\draw plot [smooth, tension=1] coordinates {(0,-3) (1.3,-1.2) (0,2)};
\node[scale=0.7] at (0,-2.3){$p_1$};
\node[scale=0.7] at (0,1.5){$p_2$};
\node[scale=0.7] at (2.5,0){1};
\node[scale=0.7] at (-2.5,0){9};
\node[scale=0.7] at (0,-4.3){$v_1$};
\node[scale=0.7] at (0,4.3){$v_6$};
\draw[dotted] (2,-4) to (0,-2);
\end{tikzpicture}
\caption{The arcs of $\mu^6(Q)$ connecting the punctures to the boundary.} \label{lotsmoremu}
\end{figure}
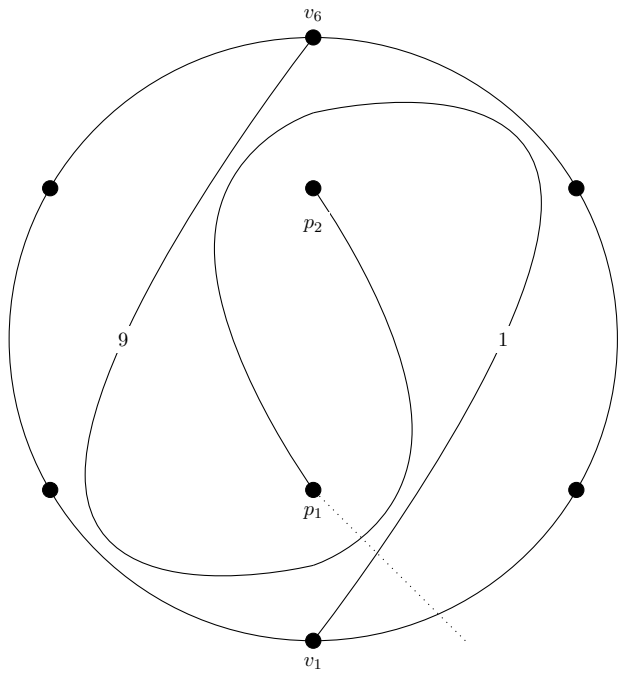

The self folded triangles are wrapping around each other more with each $\mu$. In Figure \ref{lotsmoremu} we show the arcs connected to the punctures after $6$ applications of $\mu$, which is enough to construct the rest of the triangulation, following step (ii). We see, however, that in both $Q$ and $\mu^6(Q)$ we have two self folded triangles, one at $v_1$ and $p_1$ and one at $v_6$ and $p_2$, so the description of step (i) is not sufficient. To amend this we draw a line from $p_1$ that cuts the circle somewhere between $v_1$ and $v_2$ (as shown in Figure \ref{boundarylabelling}) and consider how many times the arc from $v_1$ to $p_1$ crosses this line. This is the dashed line drawn Figure \ref{lotsmoremu}, where arc $1$ crosses it once. In this example we shall have one crossing for every $6$ applications of $\mu$. We update step (i) in our description of $\mu^l(Q)$ to capture this in the following lemma.
\begin{lemma}\label{descriptionoftriangulations}
To describe the triangulation given by applying $\mu^l$ to (\ref{Dquiver}), for any $l\in \mathbb{Z}$, we first take the boundary vertex labelling and dashed line of Figure \ref{boundarylabelling} and rotate it by  $\frac{2\pi\overline{l}}{N-2}$ clockwise, where 
\[
0\leq \overline{l}< N-2, \qquad \overline{l} \equiv l \mod N-2
\] 
while fixing $p_1$ and $p_2$. After this rotation the dashed line will still go from $p_1$ to the boundary circle between $v_1$ and $v_2$. The triangulation is then given in three steps:
\begin{enumerate}[(i)]
\item Draw a self folded triangle at $v_1$ and $p_1$, such that the arc inside the triangle crosses the dashed line $\left \lfloor{\frac{l}{N-2}}\right \rfloor$ times and travels anticlockwise (clockwise) from $v_1$ to $p_1$ if $l$ is positive (negative).
\item Draw a self-folded triangle at $v_{N-2}$ and $p_2$. There is only one choice due to step (i).
\item Draw arcs $i$ from $v_{i-2}$ to $v_{i-1}$ for $i=3,4,\ldots, N-1$. Again, there is only one choice for each of these.
\end{enumerate}
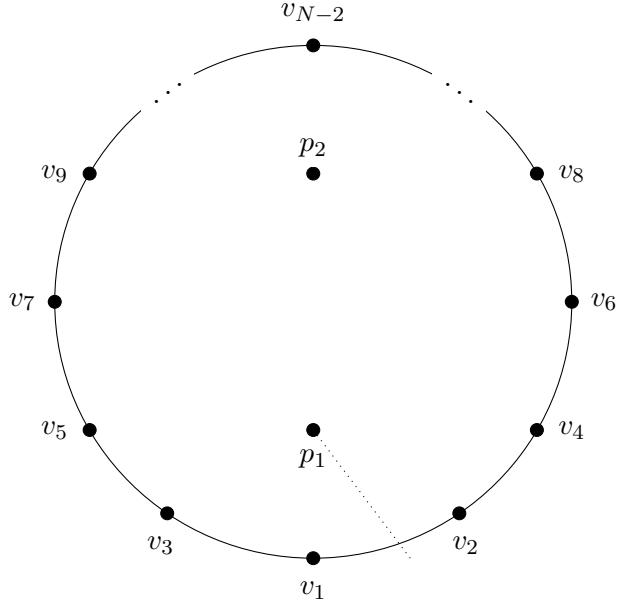
\begin{figure}
\centering
\begin{tikzpicture}[scale=0.85, every node/.style={fill=white}]
\draw (0,0) circle [radius=4.0];
\draw[fill=black] (0,4) circle [radius=0.1cm];
\draw[fill=black] (0,2) circle [radius=0.1cm];
\draw[fill=black] (0,-2) circle [radius=0.1cm];
\draw[fill=black] (0,-4) circle [radius=0.1cm];
\draw[fill=black] (3.46,2) circle [radius=0.1cm];
\draw[fill=black] (-3.46,2) circle [radius=0.1cm];
\draw[fill=black] (3.46,-2) circle [radius=0.1cm];
\draw[fill=black] (-3.46,-2) circle [radius=0.1cm];
\draw[fill=black] (4,0) circle [radius=0.1cm];
\draw[fill=black] (-4,0) circle [radius=0.1cm];
\draw[fill=black] (2.26,-3.3) circle [radius=0.1cm];
\draw[fill=black] (-2.26,-3.3) circle [radius=0.1cm];
\node[scale=1] at (-4.0,-2){$v_5$};
\node[scale=1] at (2.36,-3.8){$v_2$};
\node[scale=1] at (2.26,3.4){$\ddots$};
\node[scale=1] at (-2.26,3.4){$\iddots$};
\node[scale=1] at (-2.36,-3.8){$v_3$};
\node[scale=1] at (-4.0,2){$v_9$};
\node[scale=1] at (4.5,0){$v_6$};
\node[scale=1] at (-4.5,0){$v_7$};
\node[scale=1] at (4.0,2){$v_8$};
\node[scale=1] at (0,-4.5){$v_1$};
\node[scale=1] at (0,4.5){$v_{N-2}$};
\node[scale=1] at (4.0,-2){$v_4$};
\node[scale=1] at (0,-2.45){$p_1$};
\node[scale=1] at (0,2.4){$p_2$};
\draw[dotted] (1.5,-4) to (0,-2);
\end{tikzpicture}
\caption{Labelling around the boundary.}
\label{boundarylabelling}
\end{figure} 
\end{lemma}
Analogously to Proposition \ref{Atypebipartitebelt} we have the following result.
\begin{proposition}
The set of cluster variables obtained by the cluster map (\ref{Dtypeclustermap}), 
\[
\{X^i_n \mid i=1,2,\ldots,N+1, \quad n\in \mathbb{Z}\}
\] 
is precisely the set of cluster variables associated with $\Gamma_{1,1}\cup \Gamma_{\mathrm{punc}}$.
\begin{proof}
For this proof we take a fixed labelling of the boundary vertices as shown in Figure \ref{boundarylabelling}. We first show that every $\gamma(v_i,v_i,m)$ in $\Gamma_{1,1}$ appears due to the cluster map. From Lemma \ref{Dsurfacepossiblearcs} each of these arcs will orbit the two punctures and the line between them $m$ times before crossing the line $L$. For $i=1$ this arc appears on the outside of the self-folded triangle at $v_1$ in the quiver $\mu^{m(N-2)}(Q)$, since this encircles $L$ once for every $N-2$ applications of $\mu$. A further $\mu$ will then give $\gamma(v_3,v_3,m)$ since $\mu$ just moves the end points of $\gamma(v_1,v_1,m)$ clockwise. More applications of $\mu$ will give us $\gamma(v_i,v_i,m)$ for $i=1,2,\ldots, N-2$ (but not necessarily in this order).

To show that $\gamma(v_i,v_j,m)$ occurs for any $i,j$ we give another labelling, as shown in Figure \ref{Arooo}. We have just shown that the arc 
\[
\gamma(v_{\left \lfloor{\frac{i+j}{2}}\right \rfloor},v_{\left \lfloor{\frac{i+j}{2}}\right \rfloor},m)
\]
appears due to the cluster map. Our construction of Lemma \ref{descriptionoftriangulations} then has us draw
\[
\gamma(v_{\left \lfloor{\frac{i+j}{2}}\right \rfloor},v_{\left \lfloor{\frac{i+j}{2}}\right \rfloor+1},m)
\]
since, away from the boundary, every arc drawn in Lemma \ref{descriptionoftriangulations} is homotopic to the one drawn before. The next two we draw are
\[
\gamma(v_{\left \lfloor{\frac{i+j}{2}}\right \rfloor-1},v_{\left \lfloor{\frac{i+j}{2}}\right \rfloor+1},m), 
\qquad 
\gamma(v_{\left \lfloor{\frac{i+j}{2}}\right \rfloor-1},v_{\left \lfloor{\frac{i+j}{2}}\right \rfloor+2},m)
\]
continuing this will eventually give $\gamma(v_i,v_j,m)$. 

We also note that every arc of $\Gamma_{\mathrm{punc}}$ appears inside a $\gamma(v_i,v_i,m)$ and so will also be given by the cluster map.

Finally it is clear from Lemma \ref{descriptionoftriangulations} that every $X^i_n$ is in $\Gamma_{1,1}\cup \Gamma_{\mathrm{punc}}$.
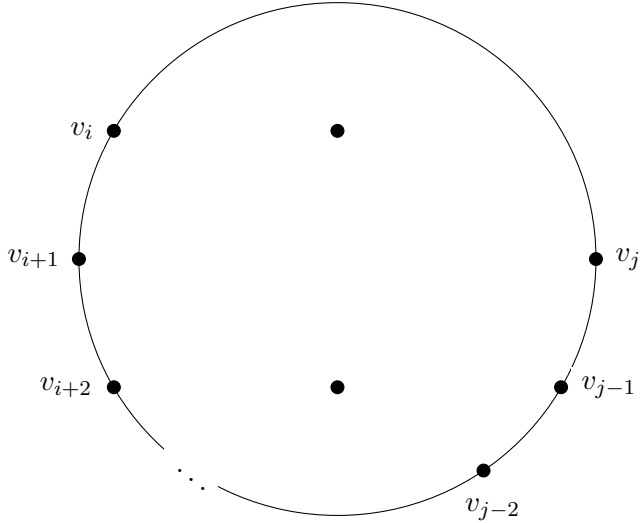
\begin{figure}
\centering
\begin{tikzpicture}[scale=0.85, every node/.style={fill=white}]
\draw (0,0) circle [radius=4.0];
\draw[fill=black] (0,2) circle [radius=0.1cm];
\draw[fill=black] (0,-2) circle [radius=0.1cm];
\draw[fill=black] (-3.46,2) circle [radius=0.1cm];
\draw[fill=black] (3.46,-2) circle [radius=0.1cm];
\draw[fill=black] (-3.46,-2) circle [radius=0.1cm];
\draw[fill=black] (4,0) circle [radius=0.1cm];
\draw[fill=black] (-4,0) circle [radius=0.1cm];
\draw[fill=black] (2.26,-3.3) circle [radius=0.1cm];
\draw[fill=black] (-2.26,-3.3) circle [radius=0.1cm];
\node[scale=1] at (-4.2,-2){$v_{i+2}$};
\node[scale=1] at (2.4,-3.9){$v_{j-2}$};
\node[scale=1] at (-2.26,-3.3){$\ddots$};
\node[scale=1] at (-3.95,2){$v_i$};
\node[scale=1] at (4.5,0){$v_{j}$};
\node[scale=1] at (-4.7,0){$v_{i+1}$};
\node[scale=1] at (4.2,-2){$v_{j-1}$};
\end{tikzpicture}
\caption{A different labelling around the boundary.}
\label{Arooo}
\end{figure} 
\end{proof}
\end{proposition} 
\begin{proposition}\label{Jproposition}
The arcs of $\Gamma_{2,0}$ live on an annulus obtained by removing a region containing the two punctures and the line between them, so we can display them on a strip, as in Figure \ref{DtypeJ}. $\Gamma_{2,0}$ is precisely the arcs ${J'}^{l}_i$ for $i=0,1,\ldots,N-3$ and $l=1,2,\ldots,N-3$. These satisfy the linear relation
\begin{equation}\label{Dtypelinear}
{J'}^{l-1}_i={J'}^{l-2}_iJ'_{i+l-2}-{J'}^{l-3}_i
\end{equation}
where $J'_{i+l-2}$ is the period $N-2$ quantity (\ref{DtypeJintro}). Hence ${J'}^l_j=D^l_1(J'_j)$, where $D^l_1(J'_j)$ is given in (\ref{determinant}), and these cluster variables form a frieze (\ref{friezeDtype}) but with final bottom row
\begin{equation*}
\begin{matrix}
\\
\ldots & D^{N-3}_1(J'_{-2}) & & D^{N-3}_1(J'_{-1}) & & D^{N-3}_1(J'_0) & & D^{N-3}_1(J'_1) & & D^{N-3}_1(J'_2) & \ldots \\
\end{matrix}
\end{equation*}
\begin{proof}
These arcs divide the disc in two such that the punctures live in the same connected component after this division. Hence the arcs avoid the part of the disc containing the line $L$ between the two punctures, so they must live in the annulus obtained by removing a small region around $L$. We draw this as in Figure \ref{DtypeJ}, where we glue along the dotted lines, and define ${J'}^{l-1}_i$ as the arc starting at $i$ and ``jumping over" $l-1$ vertices to reach $i+l$. We remark that this picture is the same as in the $\tilde{A}$ case, Figure \ref{LongerJ}, but with $0$ marked points on the internal boundary of the annulus. As such, the picture is the same as (\ref{Ptolemypic}) (with $p=1$) so we have the relation
\begin{equation}\label{Dtypelinearincomplete}
{J'}^{l-1}_j={J'}^{l-2}_j{J'}^1_{j+l-2}-{J'}^{l-3}_j.
\end{equation}
The periodic quantities for the cluster map, $J'_n$, have many equivalent expressions given in \cite{pallisterlinear}. For this proof we use 
\[
J'_n=\frac{X^3_{n+1}+X^5_n}{X^4_n}.
\]
We can obtain $J'_{-1}$ by mutating our initial quiver (triangulation), $Q$, at $3$ and then at $4$, as shown in Figure \ref{DtypeJarc}. In general $J'_n$ will be given by applying $\mu_4\mu_3$ to $\mu^{n+1}(Q)$. From this we see that $J^1_{j+l-2}=J'_{j+l-2}$ so (\ref{Dtypelinearincomplete}) becomes (\ref{Dtypelinear}). This linear relation is identical to the one in the $\tilde{A}$ case, so the determinant and frieze construction are the same as in Section \ref{Friezeconstructionsection}.
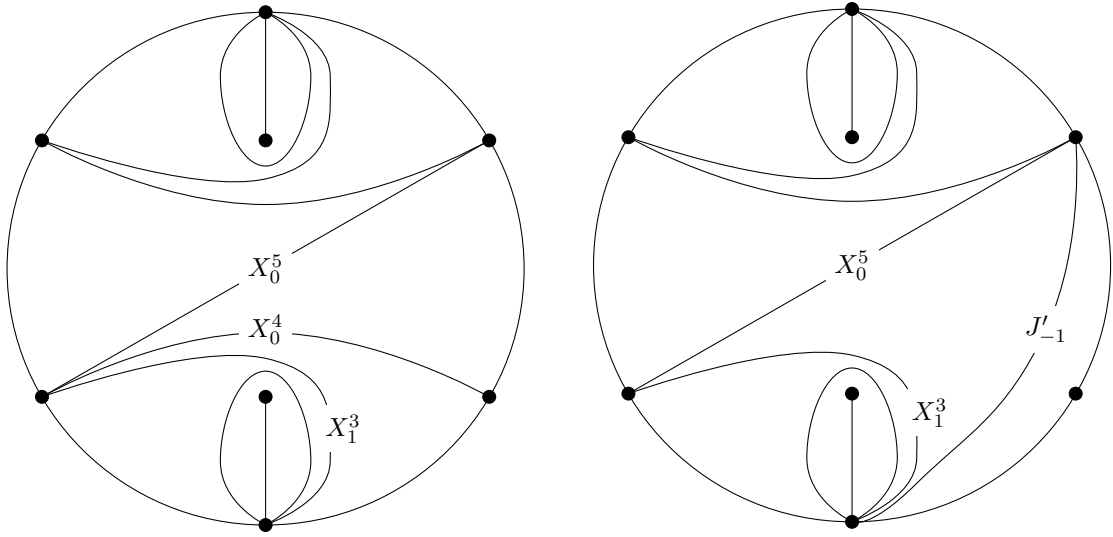
\begin{figure}
\begin{tikzpicture}[scale=0.85, every node/.style={fill=white}]
\draw (0,0) circle [radius=4.0];
\draw[fill=black] (0,4) circle [radius=0.1cm];
\draw[fill=black] (0,2) circle [radius=0.1cm];
\draw[fill=black] (0,-2) circle [radius=0.1cm];
\draw[fill=black] (0,-4) circle [radius=0.1cm];
\draw[fill=black] (3.46,2) circle [radius=0.1cm];
\draw[fill=black] (-3.46,2) circle [radius=0.1cm];
\draw[fill=black] (3.46,-2) circle [radius=0.1cm];
\draw[fill=black] (-3.46,-2) circle [radius=0.1cm];
\draw plot [smooth, tension=1] coordinates {(0,-4) (0,-2)};
\draw plot [smooth, tension=1] coordinates {(0,4) (0,2)};
\draw plot [smooth, tension=1] coordinates {(0,4) (-0.7,3) (0,1.6) (0.7,3) (0,4)};
\draw plot [smooth, tension=1] coordinates {(0,-4) (-0.7,-3) (0,-1.6) (0.7,-3) (0,-4)};
\draw plot [smooth, tension=0.9] coordinates {(0,-4) (1,-3) (0,-1.4) (-3.46,-2)};
\draw plot [smooth, tension=0.9] coordinates {(0,4) (1,3) (0,1.4) (-3.46,2)};
\draw plot [smooth, tension=0.9] coordinates {(3.46,2) (0,1) (-3.46,2)};
\draw plot [smooth, tension=0.9] coordinates {(3.46,-2) (0,-1) (-3.46,-2)};
\draw plot [smooth, tension=0.9] coordinates {(3.46,2) (0,0) (-3.46,-2)};
\node[scale=0.9] at (1.2,-2.5){$X^3_{1}$};
\node[scale=0.9] at (0,-0.95){$X^4_0$};
\node[scale=0.9] at (0,0){$X^5_0$};
\end{tikzpicture}
\qquad
\begin{tikzpicture}[scale=0.85, every node/.style={fill=white}]
\draw (0,0) circle [radius=4.0];
\draw[fill=black] (0,4) circle [radius=0.1cm];
\draw[fill=black] (0,2) circle [radius=0.1cm];
\draw[fill=black] (0,-2) circle [radius=0.1cm];
\draw[fill=black] (0,-4) circle [radius=0.1cm];
\draw[fill=black] (3.46,2) circle [radius=0.1cm];
\draw[fill=black] (-3.46,2) circle [radius=0.1cm];
\draw[fill=black] (3.46,-2) circle [radius=0.1cm];
\draw[fill=black] (-3.46,-2) circle [radius=0.1cm];
\draw plot [smooth, tension=1] coordinates {(0,-4) (0,-2)};
\draw plot [smooth, tension=1] coordinates {(0,4) (0,2)};
\draw plot [smooth, tension=1] coordinates {(0,4) (-0.7,3) (0,1.6) (0.7,3) (0,4)};
\draw plot [smooth, tension=1] coordinates {(0,-4) (-0.7,-3) (0,-1.6) (0.7,-3) (0,-4)};
\draw plot [smooth, tension=0.9] coordinates {(0,-4) (1,-3) (0,-1.4) (-3.46,-2)};
\draw plot [smooth, tension=0.9] coordinates {(0,4) (1,3) (0,1.4) (-3.46,2)};
\draw plot [smooth, tension=0.9] coordinates {(3.46,2) (0,1) (-3.46,2)};
\draw plot [smooth, tension=0.9] coordinates {(3.46,2) (3,-1) (1,-3.4) (0,-4)};
\draw plot [smooth, tension=0.9] coordinates {(3.46,2) (0,0) (-3.46,-2)};
\node[scale=0.9] at (1.2,-2.3){$X^3_{1}$};
\node[scale=0.9] at (3,-1){$J'_{-1}$};
\node[scale=0.9] at (0,0){$X^5_0$};
\end{tikzpicture}
\caption{Applying $\mu_3$ then $\mu_4$ to the initial triangulation gives the arcs $J'_{-1}$.}
\label{DtypeJarc}
\end{figure}
\begin{figure}
\centering
\begin{tikzpicture}
\draw[fill=black] (0,0) circle [radius=0.1cm];
\draw[fill=black] (2,0) circle [radius=0.1cm];
\draw[fill=black] (4,0) circle [radius=0.1cm];
\draw[fill=black] (8,0) circle [radius=0.1cm];
\draw[fill=black] (12,0) circle [radius=0.1cm];
\draw plot [smooth, tension=1] coordinates {(0,0) (4,0.7) (8,0)};
\draw plot [smooth, tension=1] coordinates {(-1,0) (2,0)};
\draw plot [smooth, tension=1] coordinates {(-1,2) (9,2)};
\draw plot [smooth, tension=1] coordinates {(11,2) (13,2)};
\draw plot [smooth, tension=1] coordinates {(11,0) (13,0)};
\draw plot [smooth, tension=1] coordinates {(2,0) (4,0)};
\draw plot [smooth, tension=1] coordinates {(7,0) (9,0)};
\node at (6,0) {$\ldots$};
\node at (2,0.9) {${J'}^{l-1}_i$};
\node at (0,-0.3) {$v_i$};
\node at (12,-0.3) {$v_i$};
\node at (2,-0.3) {$v_{i+1}$};
\node at (4,-0.3) {$v_{i+2}$};
\node at (8,-0.3) {$v_{i+l}$};
\node at (10,0) {$\ldots$};
\node at (10,2) {$\ldots$};
\draw[dotted] (-0.5,-1) to (-0.5,3);
\draw[dotted] (11.5,-1) to (11.5,3);
\end{tikzpicture}\caption{The arcs ${J'}^{l-1}_{i}$.}\label{DtypeJ}
\end{figure}
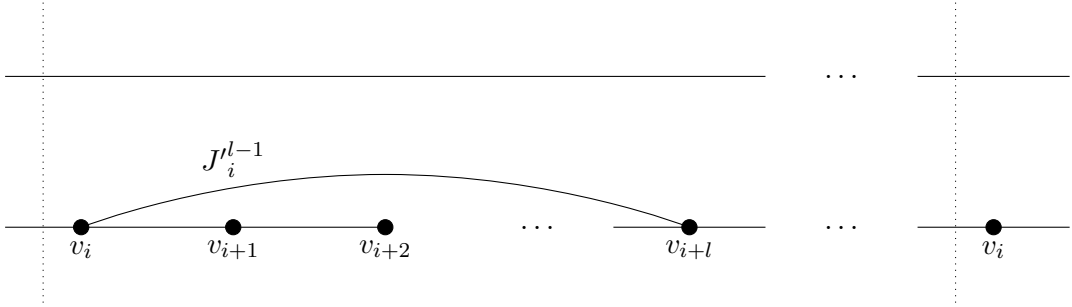 
\end{proof}
\end{proposition}
We have now found all the arcs in the $\tilde{D}$ case, except for the three only involving the punctures, $\Gamma_{\mathrm{except}}$. We collect these results in the following theorem.
\begin{theorem}
The $\tilde{D}$ type cluster variables are given by
\begin{equation}
\left\{X^i_n\:\middle|\:
\begin{aligned}
&i=1,\ldots,N+1 \\ &n\in\mathbb{Z}
\end{aligned}
\right\}
\cup
\left\{D^l_1(J'_{j})\:\middle|\:
\begin{aligned}
j=0,\ldots,N-3 \\ l=1,\ldots,N-3
\end{aligned}
\right\}
\cup 
\Gamma_{\mathrm{except}}
\end{equation} 
where the $X^i_n$ are obtained by the cluster map  (\ref{Dtypeclustermap}) and the ${J'}^{l}_i$ are defined in Proposition \ref{Jproposition}. The three arcs of $\Gamma_{\mathrm{except}}$ are shown in (\ref{3exceptionalarcs}).
\end{theorem}
We have now proven the $\tilde{D}$ part of Theorem \ref{findingallclustervarstheorem} and Proposition \ref{friezepropintro}. In the next section we complete our results for $\tilde{D}$ type with a proof of Proposition \ref{intropropclusterfrieze}.
\subsection{$\tilde{D}$ type cluster friezes}\label{Dtypeclusterfriezebit}
In $\tilde{D}$ type there are $3$ exceptional tubes, for $\lambda=1,0,\infty$. whose periods are $N-2,2,2$ respectively. We have not constructed friezes of width two, only one of width $N-2$, so our goal is to prove that our friezes agree with the friezes of \cite{assemdupont} away from the width two tubes. To do this we prove that the relation (\ref{clusterfrieze}), holds for $\lambda=1$

Firstly we follow \cite{assemdupont} to find the cluster variables associated with $B_1$ and $B'_1$. The vertex $e$ is required to be a sink, so we first mutate Figure \ref{Dquiver} at $1$ and $2$, so the quiver near $e=2$ is
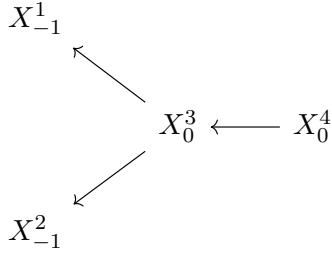
\begin{equation*}
\begin{tikzcd}
X_{-1}^1 \\
& X_0^3\arrow[ul]\arrow[dl] & X_0^4\arrow[l] \\
X_{-1}^2
\end{tikzcd}
\end{equation*}
We have $B_1=P_3/P_1$ and $B'_1=P_1[1]$ with $X_{P_1[1]}=X^1_{-1}$ and $X_{S_2}=X^2_0$. Performing $\mu_3\mu_2$ gives
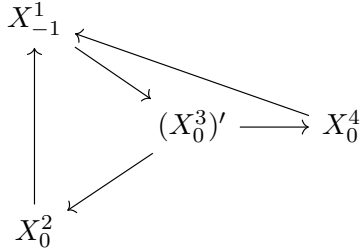
\begin{equation*}
\begin{tikzcd}
X_{-1}^1\arrow[dr] \\
& (X_0^3)'\arrow[dl]\arrow[r] & X_0^4\arrow[ull] \\
X_{0}^2\arrow[uu]
\end{tikzcd}
\end{equation*}
where $(X_0^3)'=\frac{X^2_0X^4_0+X^1_{-1}}{X^3_0}=\frac{1+X^3_0+X^2_{-1}X^1_{-1}}{X^2_{-1}X^3_0}$ is the cluster variable associated with $B_1$. We can express this as
\[
(X_0^3)'=\frac{X^1_0X^2_0X^4_0+X^1_{-1}X^1_0}{X^1_0X^3_0}=\frac{X^3_1X^3_0-1+X^3_0+1}{X^1_0X^3_0}=\frac{X^3_1+1}{X^1_0}=X^1_1
\]
so we have
\[
X_{N_{0}}=\frac{X_{B_1}+X_{B'_1}}{X_{S_1}}=\frac{X^1_1+X^1_{-1}}{X^2_0}=J'_0
\]
which gives
\[
X_{N_{0}[j]}=\frac{X_{B_1[j]}+X_{B'_1[j]}}{X_{S_1[j]}}=\frac{X^1_{1-j}+X^1_{-1-j}}{X^2_{-j}}=J'_{-j}
\]
as desired.
\clearpage
\bibliographystyle{plain}
\bibliography{Affineclusteralgebras3ways}
\end{document}